\documentclass[10pt]{article}

\usepackage{a4wide}
\usepackage{amssymb}
\usepackage{amsfonts}
\usepackage{amsmath}
\usepackage[all]{xy}
\date{}

\newtheorem{proposition}{Proposition}[section]
\newtheorem{theorem}[proposition]{Theorem}
\newtheorem{lemma}[proposition]{Lemma}

\newtheorem{corollary}[proposition]{Corollary}

\def\der{\partial }

\def\nFM0{{\nu }_{F,M_0}}
\def\nFN0{{\nu }_{F,N_0}}
\def\nGN0{{\nu }_{G,N_0}}

\def\N0{ {\bf N}_0 }

\def\ra{\rightarrow}

\def\Xpm{X^{\pm }}

\def\s{\sigma}

\def\l1{{\lambda}_1}

\def\a{\alpha}
\def\a0{ {\alpha }_0}
\def\a1{ {\alpha }_1}

\def\l{\lambda}

\def\nFGM0{{\nu }_{F,G,M_0}}

\def\nFN0{{\nu}_{F,N_0}}

\def\sm{{\sigma}^m}

\def\sm1{{\sigma}^{-1}}

\def\smtp1{{\sigma}^{-t+1}}

\def\S1{S^{-1}}

\def\Xpm1{X^{\pm 1}_1}

\def\sPM1{{\sigma }^{\pm 1}}
\def\sMP1{{\sigma }^{\mp 1 }}

\def\di{{\rm d.ind}}

\def\L{\Lambda}

\def\Ytm1{Y^{t-1}}
\def\Yim1{Y^{i-1}}

\def\CL{{\cal L}}

\def\CN{{\cal N}}
\def\CS{{\cal S}}
\def\CF{{\cal F}}

\def\ass{{\rm ass}}

\def\Aut{{\rm Aut}}

\def\ker{ {\rm ker } }

\def\SL2Z{ {\rm SL}_2({\bf Z}) }

\def\CR{ {\cal R}}

\def\th{ \theta }

\def\CL{{\cal L}}

\def\Gp1{ G^{1 , 1 } }
\def\P11{ P^{-1 , 1 } }
\def\Pp1{ P^{1 , 1 } }

\def\th{\theta}

\def\nCLsr{{}^\nu\kern-2pt {\cal L}^{\sigma , \rho  }}
\def\nP{{}^\nu \kern-2pt P}
\def\nL{{}^\nu\kern-2pt L}
\def\nLL{{}^\nu\kern-2pt \Lambda}
\def\nPsr{{}^\nu\kern-2pt P^{\sigma , \rho  }}
\def\nLsr{{}^\nu\kern-2pt L^{\sigma , \rho  }}
\def\nuCL{{}^\nu\kern-2pt  {\cal L}}
\def\nCLsr{{}^\nu\kern-2pt {\cal L}^{\sigma , \rho  }}
\def\nCL1m{{}^\nu\kern-2pt {\cal L}^{-1 , 1  }}
\def\x1nu{x^\frac{1}{\nu}}
\def\xm1nu{x^{-\frac{1}{\nu}}}

\def\rad{{\rm rad}}

\def\ob{\overline{b}}

\def\CR{ {\cal R}}

\def\CN{{\cal N}}
\def\ra{\rightarrow }

\def\CB{{\cal B}}

\def\CI{{\cal I}}

\def\CT{{\cal T}}

\def\CC{ {\cal C}}

\def\nAM0{{\nu }_{{\cal A},M_0}}
\def\nAN0{{\nu }_{{\cal A},N_0}}

\def\CR{ {\cal R }}

\def\bR{\overline{R}}

\def\ga{\mathfrak{a}}
\def\gb{\mathfrak{b}}
\def\gc{\mathfrak{c}}

\def\SL{{\rm SL}}

\def\di!{\frac{\der^i}{i!}}
\def\dik!{\frac{\der^k_i}{k!}}
\def\N{\mathbb{N}}

\def\0{\overline{0}}
\def\1{\overline{1}}

\def\Ln1{\L_{n,\overline{1}}}

\def\a1{a_{\overline{1}}}

\def\S{\Sigma}

\def\vn1{\overrightarrow{n-1}}

\def\mJ{\mathbb{J}}
\def\mI{\mathbb{I}}

\def\K1{{\rm K}_1}

\def\hmI1{\widehat{\mI_1}}
\def\tmI1{\widetilde{\mI_1}}
\def\tmJ1{\widetilde{\mJ_1}}
\def\hB1{\widehat{B_1}}
\def\hCB1{\widehat{\CB_1}}

\def\Den{{\rm Den}}

\def\Ore{{\rm Ore}}

\def\Den{{\rm Den}}

\def\Ass{{\rm Ass}}

\def\maxDen{{\rm max.Den}}

\def\llrad{{\rm l.lrad}}

\def\br{\overline{r}}

\def\bt{\overline{t}}
\def\ga{\mathfrak{a}}

\def\gll{\mathfrak{l}}

\def\CLsl{{\cal L}^s_l}
\def\Densl{{\rm Den}^s_l}

\def\glsR{\mathfrak{l}^s_R}
\def\glsRt{\mathfrak{l}^s_{R, t}}
\def\pCCR{{}'{\cal C}_R}
\def\pCCwR{{}'{\cal C}^w_R}
\def\pCCwRt{{}'{\cal C}^w_{R,t}}
\def\maxOre{{\rm max.Ore}}

\begin{document}

\author{V. V. \  Bavula %(stronglquot.tex)
}

\title{The largest strong left quotient ring of a ring}

\maketitle
%\date{}

\begin{abstract}
For an arbitrary ring $R$, the {\em largest strong left quotient ring} $Q_l^s(R)$ of $R$ and the {\em strong left localization radical} $\glsR$ are  introduced and their  properties are studied in detail. In particular, it is proved that $Q_l^s(Q_l^s(R))\simeq Q_l^s(R)$, $\gll^s_{R/\glsR}=0$ and a criterion is given for  the ring $ Q_l^s(R)$ to be  a semisimple ring. There is a canonical homomorphism from the classical left quotient ring $Q_{l, cl}(R)$ to $Q_l^s(R)$ which is not an isomorphism, in general. The objects $Q_l^s(R)$ and $\gll^s_R$ are explicitly described for several large classes of rings (semiprime left Goldie ring, left Artinian rings, rings with left Artinian left quotient ring, etc).

$\noindent $



$${\bf Contents}$$
\begin{enumerate}
\item Introduction.
\item Preliminaries, the largest strong left denominator set $T_l(R)$ of $R$ and its characterizations.
\item The  largest strong left quotient ring of a ring and its properties.
 \item The largest strong  quotient ring of a ring.
 \item The largest strong  quotient ring of $Q_l(R)$.
 \item Examples.
\end{enumerate}
\end{abstract}

%%%%%%%%%%%%%%%%%% SECTION 1 %%%%%%%%%%%%%%%%%%%%%%%%

\section{Introduction}

The aim of the paper is, for an arbitrary ring $R$, to introduce new concepts: {\em the largest strong left denominator set $T_l(R)$ of $R$, the largest strong left quotient ring $Q_l^s(R):= T_l(R)^{-1}R$ of $R$ and the strong left localization radical $\glsR $ of $R$},  and to study their properties.

In this paper, %module means a left module, and
the following notation is fixed:
\begin{itemize}
\item  $R$ is a ring with 1 and $R^*$ is its group of units;  %$\gn =\gn_R$ is its prime radical and $\Min (R)$ is the set of minimal primes of $R$;
\item   $\CC = \CC_R$  is the set of {\em regular} elements of the ring $R$ (i.e.\ $\CC$ is the set of non-zero-divisors of the ring $R$);
    \item   $\pCCR $  is the set of {\em left  regular} elements of the ring $R$, i.e.\ $\pCCR :=\{ c\in R\, | \, \ker (\cdot c)=0\}$ where $\cdot c: R\ra R$, $r\mapsto rc$;
\item   $Q=Q_{l,cl}(R):= \CC^{-1}R$ is the {\em left quotient ring}  (the {\em classical left ring of fractions}) of the ring $R$ (if it exists, i.e.\ if $\CC$ is a left Ore set) and $Q^*$ is the group of units of $Q$;
%\item   $\gn =\gn_R$ is the  prime radical of $R$ and $\nu\in \N \cup \{ \infty \}$ is its {\em nilpotency degree} ($\gn^\nu \neq 0$ but $\gn^{\nu +1}=0$);
%\item   $\bR := R/ \gn$ and $\pi: R\ra \bR$, $r\mapsto \br =r+\gn$;
%\item   $\OCC := \CC_{\bR}$ is the set of regular elements of the ring $\bR$ and $\bQ := \OCC^{-1}\bR$ is its left quotient ring;
%\item   $\CC':= \pi^{-1}(\OCC):=\{ c\in R\, \, | \, c+\gn \in \OCC\}$ and $Q':=\CC'^{-1}R$ (if it exists),

\item $\Ore_l(R):=\{ S\, | \, S$ is a left Ore set in $R\}$; \item
$\Den_l(R):=\{ S\, | \, S$ is a left denominator set in $R\}$;
\item
$\Ass_l(R):= \{ \ass (S)\, | \, S\in \Den_l(R)\}$ where $\ass
(S):= \{ r\in R \, | \, sr=0$ for some $s=s(r)\in S\}$;
\item $\Den_l(R, \ga )$ is the set of left denominator sets $S$ of $R$ with $\ass (S)=\ga$ where $\ga$ is an ideal of $R$;
     %and $\ass (S):= \{r\in R\, | \, sr=0$ for some $s\in S\}$,
        \item $S_\ga=S_\ga (R)=S_{l,\ga }(R)$
 is the {\em largest element} of the poset $(\Den_l(R, \ga ),
\subseteq )$ and $Q_\ga (R):=Q_{l,\ga }(R):=S_\ga^{-1} R$ is  the
{\em largest left quotient ring associated with} $\ga$. The fact that $S_\ga $
exists is proven in {\cite[Theorem 2.1]{larglquot}} (but also see Lemma \ref{bb23Sep13} below for the easy proof in other contexts);
\item In particular, $S_0=S_0(R)=S_{l,0}(R)$ is the largest
element of the poset $(\Den_l(R, 0), \subseteq )$, i.e.\ the {\em largest regular  left Ore set} of $R$,  and
$Q_l(R):=S_0^{-1}R$ is the {\em largest left quotient ring} of $R$ \cite{larglquot};
\item $\maxDen_l(R)$ is the set of maximal left denominator sets of $R$ (it is always a {\em non-empty} set, see \cite{larglquot}, or  Lemma \ref{bb23Sep13} below for the proof).
\end{itemize}

{\bf The largest strong left quotient ring of a ring}.
Consider the following  subsets of a ring $R$:
 The sets
\begin{eqnarray*}
\CLsl (R)&:=& \bigcap_{S\in \maxDen_l(R)}S\; \stackrel{{\rm Prop.}\, \ref{a19Sep13}.(1)}{=} \{ c\in R\, | \, \frac{c}{1}\in (S^{-1}R)^*\;\; {\rm for \; all}\;\; S\in \maxDen_l(R)\},\\
\CC^w_R&:=& \{ c\in R\, | \, \frac{c}{1}\in \CC_{S^{-1}R}\;\; {\rm for \; all}\;\; S\in \maxDen_l(R)\},   \\
 \pCCwR &:=& \{ c\in R\, | \, \frac{c}{1}\in {}'\CC_{S^{-1}R}\;\; {\rm for \; all}\;\; S\in \maxDen_l(R)\},
\end{eqnarray*}
are called respectively  {\em the set of strongly left localizable elements,   the set of  weak regular elements }  and {\em the set of  weak left regular elements} of $R$.

\begin{itemize}
\item {\rm (Proposition \ref{b19Sep13} and Proposition \ref{b22Sep13})} {\em Each of the sets $\CLsl (R)$, $\CC^w_R$ and $\pCCwR$ contains a unique largest left denominator set,  and all three largest left denominator sets coincide and are denoted by $T_l(R)$.}
\end{itemize}
The set $T_l(R)$ is called the {\em largest strong left denominator set} of $R$ and the ring $Q_l^s(R):= T_l(R)^{-1}R$ is called the {\em largest strong left quotient ring} of $R$. The ideal of $R$ given by  $\glsR :=\ass (T_l(R)):=\{ r\in R\, | \, tr=0$ for some $t\in T_l(R)\}$ is called the {\em strong left localization radical} of $R$.
 %For each ideal $\ga \in \Ass_l(R)$, there exists the {\em largest} element $S_{l, \ga}(R)$  in the poset $(\Den_l(R, \ga ) , \subseteq )$, \cite{larglquot}. The ring $Q_{l,\ga}(R):=S_{l, \ga}^{-1}R$ is called the {\em largest left quotient ring} of $R$ associated with $\ga$, \cite{larglquot}.
 In the above definitions, the adjective `strong' reflects their connections with the set $\CLsl (R)$ of strongly left localizable elements of $R$. The set $T_l(R)$  is the largest  left denominator set of $R$ that consists of elements that are invertible in {\em all} maximal left localizations of the ring $R$.

In general, for a ring $R$, its left (right; two sided) localizations, especially maximal ones, are unrelated. The intuition behind the construction of the largest strong left quotient ring of $R$ is to have the {\em largest} possible left localization of $R$  that is {\em related} to all maximal left localizations of the ring $R$, i.e. there exists a ring $R$-homomorphism (necessarily, unique) from $Q_l^s(R)$ to $S^{-1}M$ for every $S\in \maxDen_l(R)$.

 In Section \ref{EXAMPQ}, the triple
$T_l(R)$, $\gll_R^s$, $Q_l^s(R)$ is found explicitly for the following four  classes of rings:
 semiprime left Goldie rings (Theorem \ref{A29Aug14});  rings of $n\times n$ lower/upper triangular matrices with coefficients in a left Goldie domain (Theorem \ref{B29Aug14} and Theorem \ref{C29Aug14}); left Artinian rings (Theorem \ref{D29Aug14}), and  rings with left Artinian left quotient ring  (Theorem \ref{F29Aug14}). In particular,  for semiprime left Goldie rings $R$: $T_l(R)=\CC_R$, $\gll^s_R=0$ and  $Q_l^s(R)= Q_{l,cl}(R)$  (Theorem \ref{A29Aug14}). In general,  none of the three  equalities holds for the remaining three (just mentioned)  classes of rings but the results are natural and beautiful (very symmetrical), eg, for a ring $A$ such that $Q_{l,cl}(R)$ is a left Artinian ring (Theorem \ref{F29Aug14}):
 $$ T_l(A) = \bigcap_{S'\in \maxDen_l(A)}\, S'\;\; {\rm and}\;\; Q_l^s(A) \simeq \prod_{S'\in \maxDen_l(A)}S'^{-1}A.$$
 In particular, the ring $Q_l^s(A)$ has a more regular/simple structure than the ring $Q_{l,cl}(A)$.

It would be interesting to find the ring $Q_l^s(R)$ for other classes of rings. Theorem \ref{E29Aug14}.(3), which states that $Q_l^s(R) \simeq Q_l^s(Q_l(R))$ for an arbitrary ring $R$, opens a way for tackling more challenging types of rings.

 The main results of the paper are the following six theorems. The first one describes  $T_l(R)$, $Q_l^s(R)$ and $Q_l^s(R)^*$.

\begin{theorem}\label{21Sep13}%\marginpar{21Sep13}
Let $R$ be a ring, $\pi : R\ra R/ \gll^s_R$, $r\mapsto \br = r+\gll^s_R$; $\s : R\ra Q_l^s(R)$, $r\mapsto \frac{r}{1}$,  and $Q_l^s(R)^*$ be the group of units of the ring $Q_l^s(R)$. Then
\begin{enumerate}
\item $T_l(R) = S_{l,\glsR }(R)$.
\item $Q_l^s(R) = Q_{l, \gll^s_R}(R) \simeq  Q_l(R/ \gll^s_R)$.
\item $T_l(R) = \s^{-1}(Q_l^s(R)^*)$.
\item $T_l(R) = \pi^{-1} (S_{l,0}(R/\gll^s_R))$.
\item $Q_l^s(R)^*= \{ s^{-1}t\, | \, s,t \in T_l(R)\}$.
\end{enumerate}
\end{theorem}
The second one describes the objects $Q_l^s(Q_l^s(R))$,  $\gll^s_{R/ \gll^s_R}$, $T_l(R/ \gll^s_R )$, $Q_l^s(R/ \gll^s_R )$ and their connections with their counterparts for the ring $R$.

\begin{theorem}\label{B21Sep13}%\marginpar{B21Sep13}
We keep the notation of Theorem \ref{21Sep13}. Then
\begin{enumerate}
\item  $Q_l^s(Q_l^s(R))=Q_l^s(R)$.
\item $T_l(R/ \gll^s_R) = \pi (T_l(R))$ and $ T_l(R) = \pi^{-1} (T_l ( R/ \gll^s_R))$.
\item $T_l(R/ \gll^s_R) =S_{l,0}(R/ \gll^s_R) $.
\item $\gll^s_{R/ \gll^s_R}=0$.
\item $T_l(Q_l^s(R)) = Q_l^s(R)^*$  and $\gll^s_{Q_l^s(R)}=0$.
\item $\pi (\CL^s_l(R))= \CL^s_l(R/ \gll^s_R)$ and $\CL^s_l(R) = \pi^{-1} (\CL^s_l(R/ \gll^s_R))$.
    \item $Q_l^s(R/ \glsR ) = Q_l(R/ \glsR )$.
\end{enumerate}
\end{theorem}

{\bf Semisimplicity criterion for the ring $Q_l^s(R)$}. A ring is called a {\em left Goldie ring} if it does not contain infinite direct sums of nonzero left ideals and satisfies the ascending chain condition on left annihilators.

\begin{theorem}\label{Y20Sep13}%\marginpar{Y20Sep13}
Let $R$ be a ring. The following statements are equivalent.
\begin{enumerate}
\item $Q_l^s(R)$ is a semisimple ring.
\item $R/ \glsR $ is a semiprime left Goldie ring.
\item $Q_l(R/ \glsR )$ is a semisimple ring.
\item $Q_{l, cl}(R/ \glsR )$ is a semisimple ring.
\end{enumerate}
If one of the equivalent conditions holds then $$Q_l^s(R)\simeq  Q_l(R/ \glsR ) \simeq Q_{l,cl}(R/ \glsR ),$$  $T_l(R) = \pi^{-1} (\CC_{R/ \glsR })$ and $T_l(R/ \glsR )=\CC_{R/ \glsR}$  where $\pi : R\ra R/ \glsR $, $r\mapsto \br = r+\glsR$,  and $\CC_{R/ \glsR }$ is the set of regular elements of the ring $R/ \glsR $.
\end{theorem}

Goldie's Theorem \cite{Goldie-PLMS-1960} is a criterion for a ring to have  semisimple left quotient ring  (earlier,  criteria were given, by Goldie \cite{Goldie-PLMS-1958} and Lesieur and Croisot \cite{Lesieur-Croisot-1959}, for a ring to have  a simple Artinian left quotient  ring). Recently, the author \cite{Bav-Crit-S-Simp-lQuot} has  given several more new  criteria.
 For a left Noetherian ring which has
a left quotient ring, Talintyre \cite{Talintyre-QuotRinMin-66} has established necessary and sufficient
conditions for the left quotient ring to be left Artinian. Small \cite{Small-ArtQuotRings-66, Small-CorrArtQuotRings-66},    Robson \cite{Robson-ArtQuotRings-67}, and later  Tachikawa \cite{Tachikawa-AutQuotR-71} and Hajarnavis \cite{Hajarnavis-ThmSmall-72}, and recently the author \cite{Bav-genGoldie}  have given different criteria  for a ring to have a left Artinian left quotient ring.

{\bf Semisimplicity criterion for the ring $Q_{l,cl}(R)$}. The statement of Goldie's Theorem is a semisimplicity criterion for the ring $Q_{l,cl}(R)$ which  states that the ring $Q_{l,cl}(R)$ is a semisimple ring iff $R$ is a semiprime left Goldie ring. Recently, four new criteria for semisimplicity of $Q_{l,cl}(R)$ are given in \cite{Bav-Crit-S-Simp-lQuot} using completely different ideas and approach. Below, another semisimplicity criterion for $Q_{l,cl}(R)$ is given via $Q_l^s(R)$ and $\glsR$.

\begin{theorem}\label{27Sep13}%\marginpar{27Sep13}
Let $R$ be a ring. The following statements are equivalent.
\begin{enumerate}
\item $Q_{l, cl}(R)$ is a semisimple ring.
\item $Q_l(R)$ is a semisimple ring.
\item
\begin{enumerate}
\item $Q_l^s(R)$ is a semisimple ring.
%\item $\CC_R\subseteq T_l(R)$.
\item $\glsR =0$.
\end{enumerate}
\end{enumerate}
If one of the equivalent conditions 1--3  holds then $Q_{l,cl}(R)\simeq Q_l(R)\simeq  Q_l^s(R) $ and  $\CC_R = T_l(R)=\CLsl (R)$.
\end{theorem}

\begin{theorem}\label{A29Sep13}%\marginpar{A29Sep13}
Let $R$ be a ring. Then, for all $\ga \in \Ass_l(R)$ with $\ga \subseteq \glsR$,   $S_{l,\ga}(R)\subseteq T_l(R)$, and so there is a ring $R$-homomorphism $Q_{l , \ga}(R)\ra Q_l^s(R)$, $s^{-1}r\mapsto s^{-1}r$.
\end{theorem}

  For an arbitrary ring $R$, Theorem \ref{E29Aug14} reveals natural and tight connections between triples $T_l(R)$, $\gll_R^s$, $Q_l^s(R)$ and  $T_l(Q_l(R))$, $\gll_{Q_l(R)}^s$, $Q_l^s(Q_l(R))$.

\begin{theorem}\label{E29Aug14}%\marginpar{E29Aug14}
Let $R$ be a ring. Then
\begin{enumerate}
\item $T_l(Q_l(R)) = Q_l(R)^*\, T_l(R)= \{ s^{-1}t\, | \, s\in S_l(R), t\in T_l(R)\}$ and $T_l(R) = R\cap T_l(Q_l(R))$.
\item $\gll^s_{Q_l(R)}=S_l(R)^{-1}\gll^s_R$ and $\gll^s_R= R\cap \gll^s_{ Q_l(R)}$.
\item $Q^s_l(R)\simeq Q^s_l(Q_l(R))$.
\end{enumerate}
\end{theorem}

The paper is organized as  follows. In Section \ref{PRLM}, we prove Proposition \ref{b19Sep13} and Proposition \ref{b22Sep13} (mentioned above). We show that that  $S_{l,0}(R)\subseteq T_l(R)$ (Lemma \ref{c19Sep13}.(3)) and as a result there is a canonical homomorphism
$$ \th : Q_l(R)\ra Q_l^s(R), \;\; s^{-1}r\mapsto s^{-1}r, \;\; (s\in S_l(R), \; r\in R).$$
The lemma below is a criterion for the homomorphism $\th$ to be an isomorphism.
\begin{itemize}
\item {(\rm Lemma \ref{a15Sep13})} $S_{l,0}(R) = T_l(R)$ {\em iff $\th$ is an isomorphism iff} $\glsR =0$.
\end{itemize}
In Section \ref{TSLQRP}, proofs of Theorems \ref{21Sep13}--\ref{A29Sep13}  %Theorems \ref{B21Sep13}, Theorem \ref{Y20Sep13} and Theorem \ref{27Sep13}
 are given. In Section \ref{T4SQR}, the two-sided theory  (i.e. about left and right denominators sets)  is developed and analogous results to the  five theorems above  are proved.
In Section \ref{LSQQlR},  Theorem \ref{E29Aug14} is proved.

%$\noindent $

%%%%%%%%%%%%%%%%%% SECTION 2  %%%%%%%%%%%%%%%%%%%%%%%%

\section{Preliminaries, the largest strong left denominator set $T_l(R)$ of $R$ and its characterizations }\label{PRLM}%\marginpar{PRLM}

In this section, for reader's convenience we collect necessary results that are used in the proofs of this paper. Several characterizations (Proposition \ref{b19Sep13} and Proposition \ref{b22Sep13}) of $T_l(R)$ are given. A criterion is given for the inclusion $S_{l,0}(R)\subseteq T_l(R)$ (which always holds by Lemma \ref{c19Sep13}.(3)) to be an equality,   and for the canonical ring homomorphism $Q_{l, cl}(R)\ra Q_l^s(R)$ to be an isomorphism (Lemma \ref{a15Sep13}).

 More results on localizations of rings (and some of the missed standard definitions) the reader can find in \cite{Jategaonkar-LocNRings}, \cite{MR} and \cite{Stenstrom-RingQuot}.

%{\bf Notation}:

%\begin{itemize}
%\item $\Ore_l(R):=\{ S\, | \, S$ is a left Ore set in $R\}$; \item
%$\Den_l(R):=\{ S\, | \, S$ is a left denominator set in $R\}$;
%\item $\Loc_l(R):= \{ S^{-1}R\, | \, S\in \Den_l(R)\}$; \item
%$\Ass_l(R):= \{ \ass (S)\, | \, S\in \Den_l(R)\}$ where $\ass
%(S):= \{ r\in R \, | \, sr=0$ for some $s=s(r)\in S\}$;
%\item $\Den_l(R, \ga ) := \{ S\in \Den_l(R)\, | \, \ass (S)=\ga \}$ where $\ga \in \Ass_l(R)$;
%\item $S_\ga=S_\ga (R)=S_{l,\ga }(R)$
% is the {\em largest element} of the poset $(\Den_l(R, \ga ),
%\subseteq )$ and $Q_\ga (R):=Q_{l,\ga }(R):=S_\ga^{-1} R$ is  the
%{\em largest left quotient ring associated to} $\ga$, $S_\ga $
%exists (Theorem %\ref{3Jul10}.(2));
% 2.1, \cite{larglquot});
%\item In particular, $S_0=S_0(R)=S_{l,0}(R)$ is the largest
%element of the poset $(\Den_l(R, 0), \subseteq )$, i.e. the largest regular  left Ore set of $R$,  and
%$Q_l(R):=S_0^{-1}R$ is the {\em largest left quotient ring} of $R$ \cite{larglquot};
%\item $\Loc_l(R, \ga ):= \{ S^{-1}R\, | \, S\in \Den_l(R, \ga
%)\}$.
%\end{itemize}

$\noindent $

{\bf The largest regular left Ore set and the largest left
quotient ring of a ring}. Let $R$ be a ring. A {\em
multiplicatively closed subset} $S$ of $R$ or a {\em
 multiplicative subset} of $R$ (i.e. a multiplicative sub-semigroup of $(R,
\cdot )$ such that $1\in S$ and $0\not\in S$) is said to be a {\em
left Ore set} if it satisfies the {\em left Ore condition}: for
each $r\in R$ and
 $s\in S$, $ Sr\bigcap Rs\neq \emptyset $.
Let $\Ore_l(R)$ be the set of all left Ore sets of $R$.
  For each $S\in \Ore_l(R)$ the set $\ass (S) :=\{ r\in
R\, | \, sr=0 \;\; {\rm for\;  some}\;\; s\in S\}$  is an ideal of
the ring $R$.

%$\noindent $

A left Ore set $S$ is called a {\em left denominator set} of the
ring $R$ if $rs=0$ for some elements $ r\in R$ and $s\in S$ implies
$tr=0$ for some element $t\in S$, i.e.\ $r\in \ass (S)$. Let
$\Den_l(R)$ be the set of all left denominator sets of $R$. For
$S\in \Den_l(R)$, let $S^{-1}R=\{ s^{-1}r\, | \, s\in S, r\in R\}$
be the {\em left localization} of the ring $R$ at $S$ (the {\em
left quotient ring} of $R$ at $S$). Let us stress that in Ore's method of localization one can localize {\em precisely} at left denominator sets.

In general, the set $\CC$ of regular elements of a ring $R$ is
neither a left nor right Ore set of the ring $R$ and as a
 result neither the left nor right classical  quotient ring ($Q_{l,cl}(R):=\CC^{-1}R$ and
 $Q_{r,cl}(R):=R\CC^{-1}$) exists.
 Remarkably, there  exists a largest
 regular left Ore set $S_0= S_{l,0} = S_{l,0}(R)$, \cite{larglquot}. This means that the set $S_{l,0}(R)$ is an Ore set of
 the ring $R$ that consists
 of regular elements (i.e.\ $S_{l,0}(R)\subseteq \CC$) and contains all the left Ore sets in $R$ that consist of
 regular elements. Also, there exists a largest regular (left and right) Ore set $S_{l,r,0}(R)$ of any ring $R$.
 In general, all the sets $\CC$, $S_{l,0}(R)$, $S_{r,0}(R)$ and $S_{l,r,0}(R)$ are distinct.  For example, these sets are different
   for   the ring $\mI_1= K\langle x, \der , \int\rangle$ of polynomial integro-differential operators  over a field $K$ of characteristic zero,    \cite{Bav-intdifline}. In  \cite{Bav-intdifline},  these four sets are found explicitly for $R=\mI_1$.

$\noindent $

{\it Definition}: Following the terminology of  \cite{Bav-intdifline}, \cite{larglquot},     we call the ring
$$Q_l(R):= S_{l,0}(R)^{-1}R$$ (respectively, $Q_r(R):=RS_{r,0}(R)^{-1}$ and
$Q(R):= S_{l,r,0}(R)^{-1}R\simeq RS_{l,r,0}(R)^{-1}$)
the {\em largest left} (respectively, {\em right and two-sided})
{\em quotient ring} of the ring $R$.

$\noindent $

 In general, the rings $Q_l(R)$, $Q_r(R)$ and $Q(R)$
are not isomorphic, for example, for $R= \mI_1$ as shown in Section 8 of  \cite{Bav-intdifline}.  The next
theorem gives various properties of the ring $Q_l(R)$. In
particular, it describes its group of units.

%$\noindent $

\begin{theorem}\label{4Jul10}\cite{larglquot}%\marginpar{4Jul10}
%\cite{larglquot}
\begin{enumerate}
\item $ S_0 (Q_l(R))= Q_l(R)^*$ {\em and} $S_0(Q_l(R))\cap R=
S_0(R)$.
 \item $Q_l(R)^*= \langle S_0(R), S_0(R)^{-1}\rangle$, {\em i.e. the
 group of units of the ring $Q_l(R)$ is generated by the sets
 $S_0(R)$ and} $S_0(R)^{-1}:= \{ s^{-1} \, | \, s\in S_0(R)\}$.
 \item $Q_l(R)^* = \{ s^{-1}t\, | \, s,t\in S_0(R)\}$.
 \item $Q_l(Q_l(R))=Q_l(R)$.
\end{enumerate}
\end{theorem}

{\bf The maximal left denominator sets and the maximal left localizations  of a ring}. The set $(\Den_l(R), \subseteq )$ is a poset (partially ordered
set). In \cite{larglquot}, it is proved  that the set
$\maxDen_l(R)$ of its maximal elements is a {\em non-empty} set.

$\noindent $

{\it Definition}, \cite{larglquot}. An element $S$ of the set
$\maxDen_l(R)$ is called a {\em maximal left denominator set} of
the ring $R$ and the ring $S^{-1}R$ is called a {\em maximal left
quotient ring} of the ring $R$ or a {\em maximal left localization
ring} of the ring $R$. The intersection
%\marginpar{llradR}
\begin{equation}\label{llradR}
\gll_R:=\llrad (R) := \bigcap_{S\in \maxDen_l(R)} \ass (S)
\end{equation}
is called the {\em left localization radical } of the ring $R$,
\cite{larglquot}.

$\noindent $

% For a ring $R$, there is the canonical exact
%sequence %\marginpar{llRseq}
%\begin{equation}\label{llRseq}
%0\ra \gll_R \ra R\stackrel{\s }{\ra} \prod_{S\in \maxDen_l(R)}S^{-1}R, \;\; \s := \prod_{S\in \maxDen_l(R)}\, \s_S,
%\end{equation}
%where $\s_S:R\ra S^{-1}R$, $r\mapsto \frac{r}{1}$ .% For a ring $R$ with a semisimple left quotient ring, Theorem \ref{1Jan13} shows that the left localization radical $\gll_R$ coincides with the prime radical $\gn_R$ of $R$: $\gll_R=\bigcap_{\gp\in \Min (R)}\gp =\gn_R= 0$. In general, this  is not the case even for left Artinian rings \cite{Bav-LocArtRing}.

{\bf Properties of the maximal left quotient rings of a ring}.
The next theorem describes various properties of the maximal left
quotient rings of a ring. In particular, their groups of units and
their largest left quotient rings. %The key moment in the proof is to use Theorem \ref{4Jul10}.
 It is the key fact in the proof of the characterization of the set $ \CLsl$ (Proposition \ref{a19Sep13}).
\begin{theorem}\label{15Nov10}\cite{larglquot}%\marginpar{15Nov10}
 $\,$  Let $S\in \maxDen_l(R)$, $A= S^{-1}R$, $A^*$ be
the group of units of the ring $A$; $\ga := \ass (S)$, $\pi_\ga
:R\ra R/ \ga $, $ a\mapsto a+\ga$, and $\s_\ga : R\ra A$, $
r\mapsto \frac{r}{1}$. Then
\begin{enumerate}
\item $S=S_\ga (R)$, $S= \pi_\ga^{-1} (S_0(R/\ga ))$, $ \pi_\ga
(S) = S_0(R/ \ga )$ and $A= S_0( R/\ga )^{-1} R/ \ga = Q_l(R/ \ga
)$. \item  $S_0(A) = A^*$ and $S_0(A) \cap (R/ \ga )= S_0( R/ \ga
)$. \item $S= \s_\ga^{-1}(A^*)$. \item $A^* = \langle \pi_\ga (S)
, \pi_\ga (S)^{-1} \rangle$, i.e. the group of units of the ring
$A$ is generated by the sets $\pi_\ga (S)$ and $\pi_\ga (S)^{-1}:=
\{ \pi_\ga (s)^{-1} \, | \, s\in S\}$. \item $A^* = \{ \pi_\ga
(s)^{-1}\pi_\ga ( t) \, |\, s, t\in S\}$. \item $Q_l(A) = A$ and
$\Ass_l(A) = \{ 0\}$.     In particular, if $T\in \Den_l(A)$
then  $T\subseteq A^*$.
\end{enumerate}
\end{theorem}

Theorems \ref{4Jul10} and \ref{15Nov10} are used in many proofs in this paper.

$\noindent $

{\it Definition}, \cite{Bav-Crit-S-Simp-lQuot}. The sets
$$ \CL_l(R):= \bigcup_{S\in \maxDen_l(R)}S\;\; {\rm and}\;\; \CN\CL_l(R):=R\backslash \CL_l(R)$$ are called the sets of {\em left localizable} and {\em left non-localizable elements} of $R$, respectively, and the intersection
$$\CL^s_l (R) :=\bigcap_{S\in \maxDen_l(R)} S$$
 is called the {\em set of strongly (or completely)  left localizable elements} of $R$. Clearly, $\CLsl (R)$ is a multiplicative set and
 %\marginpar{RsSl}
\begin{equation}\label{RsSl}
 R^*\subseteq \CL^s_l(R)
\end{equation}
 since $R^*\subseteq S$ for all $S\in \maxDen_l(R)$, by Theorem \ref{15Nov10}.(3).
By Proposition \ref{A8Dec12}.(1),
%\marginpar{SlinSL}
\begin{equation}\label{SlinSL}
S_{l,0}(R)\subseteq  \CL^s_l(R).
\end{equation}
In particular, if the set $\CC_R$ of regular elements of the ring $R$ is a left Ore set then $\CC_R = S_{l,0}(R)$ and so
%\marginpar{1SlinSL}
\begin{equation}\label{1SlinSL}
\CC_R \subseteq  \CL^s_l(R).
\end{equation}
The next proposition is a characterization of the set $\CL^s_l(R)$ which says that the set $\CLsl (R)$ contains {\em precisely} the elements of the ring $R$ that are units in all maximal left localizations of $R$.

\begin{proposition}\label{a19Sep13}%\marginpar{a19Sep13}
Let $R$ be a ring. Then
\begin{enumerate}
\item $\CL^s_l(R) = \{ s\in R\,| \, \frac{s}{1}\in (S^{-1}R)^*$ for all $S\in \maxDen_l(R)\}$ where $(S^{-1}R)^*$ is the group of units of the ring $S^{-1}R$.
\item For all automorphisms $\s \in \Aut (R)$, $\s (\CLsl (R))= \CLsl (R)$.
\end{enumerate}
\end{proposition}

{\it Proof}. 1.  Let $\CR$ be the RHS of the claimed  equality. By the very definition of the set $\CL^s_l(R)$,   we have the inclusion $\CL^s_l(R) \subseteq \CR$. Conversely, let $s\in \CR$ and $\s_S :R\ra S^{-1}R$, $r\mapsto \frac{r}{1}$, where $S\in \maxDen_l(R)$. Then $s\in \s^{-1}_S((S^{-1}R)^*)= S$ for all $S\in \maxDen_l(R)$ (Theorem \ref{15Nov10} parts 2 and 3), hence $s\in \CL^s_l(R)$.

2. Obvious.   $\Box $

$\noindent $

Let $R$ be a ring. Let $S,T$ be submonoids of the multiplicative monoid $(R, \cdot )$. We denote by $ST$ the {\em submonoid}  of $(R, \cdot )$ generated by $S$ and $T$. This notation should not be confused with the product of two sets  which is {\em not} used in this paper. The next result is a criterion for the set $ST$ to be a left Ore (denominator) set.

\begin{lemma}\label{a14Sep13}%\marginpar{a14Sep13}

\begin{enumerate}
\item Let $S,T\in \Ore_l(R)$. If $0\not\in ST$ then $ST \in \Ore_l (R)$.  \item Let $S,T\in \Den_l(R)$. If $0\not\in ST$ then $ST \in \Den_l (R)$.
    \item Statements 1 and 2 hold also for Ore sets  and denominator sets, respectively.
\end{enumerate}
\end{lemma}

{\it Proof}. 1. Since $0\not\in ST$, the set $P:=ST$ is multiplicative. It remains to show that the left Ore condition holds for $P$.  Given an element $p=s_1t_1\cdots s_nt_n\in P$ and $r\in R$ (where $s_i\in S$ and $t_i\in T$) we have to find elements $p'\in P$ and $r'\in R$ such that $p'r=r'p$. There are  elements $t_n'\in T$ and $r_n'\in R$ such that $t_n'r=r_n't_n$. Similarly, $s_n'r_n' = r_n''s_n$ for some $s_n'\in S$ and $r_n''\in R$. Hence, $s_n't_n'r= r_n''s_nt_n$. Then repeating these two steps $n-1$ more times we find elements $s_i'\in S$, $t_i'\in T$ and $ r'\in R$ such that
$$ s_1't_1'\cdots s_n't_n'r=r's_1t_1\cdots s_nt_n.$$
So,  it suffices to take $p' = s_1't_1'\cdots s_n't_n'$.

2. By statement 1, it remains  to show that if $rp =0$ for some elements $ r\in R$ and $p= s_1t_1\cdots s_nt_n\in P$ then $p'r=0$ for some $p'\in P$. $0=rp= (r s_1t_1\cdots s_n)t_n \Rightarrow t_n'r s_1t_1\cdots s_n=0$ for some element $t_n'\in T$. Similarly, $s_n' t_n'r s_1t_1\cdots s_{n-1}t_{n-1}=0$ for some element $s_n'\in S$. Repeating the same two steps $n-1$ more times we have $
s_1't_1'\cdots s_n't_n'r=0$ for some elements $s_i'\in S$ and $t_i'\in T$. It suffices to take $p'=s_1't_1'\cdots s_n't_n'$.

3. Statement 3 follows from statements 1 and 2. $\Box $

$\noindent $

{\bf Criterion for a left Ore/denominator  set to be maximal}. There are posets $(\Ore_l(R), \subseteq )$, $(\Den_l(R), \subseteq )$, $(\Ore (R), \subseteq )$ and  $(\Den (R), \subseteq )$. The next lemma states that the sets of maximal elements of these posets are non-empty sets.

\begin{lemma}\label{bb23Sep13}%\marginpar{bb23Sep13}
Let $R$ be a ring.
\begin{enumerate}
\item The set $\maxOre_l(R)$ of maximal left Ore sets in $R$ is a non-empty set.
\item The set $\maxDen_l (R)$ of  maximal left  denominator  sets  in $R$ is a non-empty set.
\item The set $\maxOre (R)$ of maximal (left and right) Ore sets  in $R$ is a non-empty set.
    \item The set $\maxDen (R)$ of maximal (left and right) denominator  sets  in $R$ is a non-empty set.
    \end{enumerate}
\end{lemma}

{\it Proof}. All statements follow at once from Zorn's Lemma and the fact that given a linearly ordered chain of left (resp. left and right) Ore sets [resp. denominator sets] then their union  is a  left (resp. left and right) Ore set  [resp. a denominator set].  $\Box $

$\noindent $

The next proposition is a criterion for a left Ore/denominator set to be a {\em maximal} left Ore/denominator set.

\begin{proposition}\label{b23Sep13}%\marginpar{b23Sep13}
Let $R$ be a ring.
\begin{enumerate}
\item Let $S\in \Ore_l(R)$ (resp. $S\in \Ore (R)$). Then $S\in \maxOre_l(R)$ (resp. $S\in {\rm max.Ore} (R)$) iff $0\in ST$ for all $T\in \Ore_l(R)$ such that $T\not\subseteq S$ (resp. $T\in \Ore (R)$).

\item Let $S\in \Den_l(R)$ (resp. $S\in \Den (R)$). Then $S\in \maxDen_l(R)$ (resp. $S\in \maxDen (R)$) iff $0\in ST$ for all $T\in \Den_l(R)$ such that $T\not\subseteq S$ (resp. $T\in \Den (R)$).
    \end{enumerate}
\end{proposition}

{\it Proof}. 1. Statement 1 follows from Lemma \ref{a14Sep13}.(1,3) and the inclusion $S\subseteq ST$.

2. Statement 1 follows from Lemma \ref{a14Sep13}.(2,3) and the inclusion $S\subseteq ST$.  $\Box $

$\noindent $

Let $\{S_i\, | \, i\in I\}\subseteq \Ore_l(R)$, $I\neq \emptyset$, $\CF := \{ J\subseteq I\, | \, 1\leq |J|<\infty\}$ and
%\marginpar{VSi}
\begin{equation}\label{VSi}
\bigvee_{i\in I} S_i :=\bigcup_{F\in \CF}\prod_{i\in F}S_i.
\end{equation}

\begin{lemma}\label{a18Sep13}%\marginpar{a18Sep13}

\begin{enumerate}
\item Let $\{S_i\, | \, i\in I\}\subseteq \Ore_l(R)$. Suppose that $0\not\in \prod_{i\in F}S_i$ for all non-empty finite subsets $F\subseteq I$.  Then  $\bigvee_{i\in I} S_i$ is the  least upper bound of $\{S_i\, | \, i\in I\}$ in $\Ore_l(R)$.
\item Let $\{S_i\, | \, i\in I\}\subseteq \Den_l(R)$. Suppose that $0\not\in \prod_{i\in F}S_i$ for all non-empty finite subsets $F\subseteq I$.  Then  $\bigvee_{i\in I} S_i$ is the  least upper bound of $\{S_i\, | \, i\in I\}$ in $\Den_l(R)$.
    \item Statements 1 and 2 hold also for Ore sets and denominator sets,  respectively.
\end{enumerate}
\end{lemma}

{\it Proof}. 1.   By Lemma \ref{a14Sep13}.(1), $\bigvee_{i\in I} S_i\in \Ore_l(R)$. Now, statement 1 is obvious.

2.   By Lemma \ref{a14Sep13}.(2), $\bigvee_{i\in I} S_i\in \Den_l(R)$. Now, statement 2 is obvious.

3. Statement 3 follows from statements 1 and 2.  $\Box $

$\noindent $

{\bf The largest strong left denominator set $T_l(R)$}.
 The set $\Densl (R):= \{ T\in \Den_l(R)\, | \, T\subseteq \CL^s_l(R)\}$ is a {\em non-empty} set since $ R^*\in \Densl (R)$. The elements of $\Densl (R)$ are called the {\em strong left denominator sets} of $R$ and the rings $T^{-1}R$ where $T\in \Densl (R)$ are called the {\em strong left quotient rings} or the {\em  strong left localizations} of $R$. The next proposition  shows that the set of maximal elements (w.r.t. to inclusion) ${\rm max.Den}^s_l(R)$ is a non-empty set. Moreover, it contains a single element. Namely,

 %\marginpar{TlR}
\begin{equation}\label{TlR}
T_l(R):=\bigcup_{S\in \Densl (R)}S.
\end{equation}

\begin{proposition}\label{b19Sep13}%\marginpar{b19Sep13}
Let $R$ be a ring. Then ${\rm max.Den}^s_l(R)=\{T_l(R)\}$.
\end{proposition}

{\it Proof}. (i)  {\em For all $S,T\in \Densl (R)$, $0\not\in ST$ where $ST$ is the semigroup of $(R, \cdot )$ generated by the sets $S$ and} $T$: Suppose that $0\in ST$ for some $S,T\in \Densl (R)$, i.e.\ $s_1t_1s_2t_2\cdots s_nt_n=0$ for some elements $s_i\in  S$ and $t_i\in T$. Take $\CS\in \maxDen_l(R)$. Then  $S,T\subseteq \CS$ and so $ 0\neq s_1t_1s_2t_2\cdots s_nt_n/1\in \CS^{-1}R$, a contradiction.

(ii) $\bigvee_{S\in \Densl (R)}S= T_l(R)$ (see (\ref{VSi})): This follows from the fact that $\Densl (R)$ is a monoid, by (i).

(iii) ${\rm max.Den}^s_l(R) = \{ T_l(R)\}$: By (ii) and Lemma  \ref{a18Sep13}.(2), $T_l(R)$ is the least  upper bound of $\Densl (R)$
 in the set $\Den_l (R)$. Since  $T_l(R)\in \Densl (R)$,  $T_l(R)$ is the largest element in $\Densl (R)$. $\Box$

$\noindent $

So, $T_l(R)$ is the {\em largest strong denominator set} of $R$.

 $\noindent $

{\it Definition}. The ideal of the ring $R$,
$$ \gll^s_R := \ass (T_l(R))= \bigcup_{S\in \Densl (R)}\ass (S)$$ is called the {\em strong left localization radical } of the ring $R$.

$\noindent $

The set $\Densl (R)$ is invariant under the action of the group of automorphisms $\Aut (R)$ of the ring $R$.
\begin{lemma}\label{bb19Sep13}%\marginpar{bb19Sep13}

\begin{enumerate}
\item For all automorphisms $\s \in \Aut (R)$, $\s (T_l(R))= T_l(R)$ and  $ \s (\gll^s_R)= \gll^s_R$.
\item $\gll^s_R\subseteq \gll_R$ where $ \gll_R :=\bigcap_{S\in \maxDen_l(R)}\ass (S)$.
\end{enumerate}
\end{lemma}

{\it Proof}. 1. Obvious.

2. For all $S\in \maxDen_l(R)$, $T_l(R)\subseteq S$, and so $\ass (T_l(R))\subseteq \ass (S)$. Therefore,  $\gll^s_R\subseteq  \gll_R$. $\Box $

%$\noindent $

The next result shows that the sets $T_l(R)$ and $ \CL_l^s(R)$ are closed under addition of elements of $\gll^s_R$ and $\gll_R$, respectively.

\begin{lemma}\label{c19Sep13}%\marginpar{c19Sep13}
Let $R$ be a ring. Then
\begin{enumerate}
\item $T_l(R)+\gll^s_R =T_l(R)$.
\item  $\CLsl (R) + \gll_R\subseteq \CLsl (R)$. In particular, $T_l(R) + \gll_R\subseteq \CL^s_l(R)$.
\item $R^*\subseteq S_{l,0}(R)\subseteq T_l(R)$.
\end{enumerate}
\end{lemma}

{\it Proof}. 1. Let $T_l := T_l(R)$. By Corollary \ref{b21Sep13}.(1) (or Lemma \ref{a21Sep13}) (see below), $T:= T_l+\gll^s_R\in \Den_l(R)$. To finish the  proof of statement 1 it suffices to show that $T\subseteq \CL^s_l(R)$ (since then $T\subseteq T_l$, as $T_l$ is the largest element in $\CL^s_l(R)$, Proposition \ref{b19Sep13}). Let $t\in T_l$ and $ a\in \gll^s_R$. We have to show that $t+a\in \CL^s_l(R)$. For all $S\in \maxDen_l(R)$, $T_l\subseteq S$ and so $\gll^s_R\subseteq \ass (S)$. Hence,
$$ S^{-1}R \ni \frac{t+a}{1}=\frac{t}{1}\in (S^{-1}R)^*.$$
By Proposition \ref{a19Sep13}.(1), $t+a\in \CL^s_l(R)$.

2. Let $t\in \CLsl (R)$ and $l\in \gll_R$. For all $S\in \maxDen_l(R)$, $S^{-1}R \ni \frac{t+l}{1}=\frac{t}{1}\in (S^{-1}R)^*$. By  Proposition \ref{a19Sep13}, $t+l\in \CL^s_l(R)$.

3. By the inclusion (\ref{SlinSL}), $S_{l,0}(R)\in \Densl (R)$. Hence, $ S_{l,0}(R) \subseteq T_l(R)$, by Proposition \ref{b19Sep13}. By Theorem \ref{4Jul10}.(1), $R^*\subseteq S_{l,0}(R)$.  $\Box $

%$\noindent $

The next lemma  shows that under natural conditions pre-images of left denominator sets are also left denominator sets.

\begin{lemma}\label{a21Sep13}%\marginpar{a21Sep13}
Let $R$ be a ring, $S\in \Den_l(R, \ga )$, $\pi : R\ra R/ \ga$, $r\mapsto \br := r+\ga$,  $T\in \Den_l(R/ \ga , \gb / \ga )$ where $\gb$ is an ideal of $R$ such $\ga \subseteq \gb$, and $ T':= \pi^{-1}(T)$. If $S\subseteq T'$ then $ T'\in \Den_l(R, \gb )$.
\end{lemma}

{\it Proof}. (i) $T'\in \Ore_l(R)$: Given elements $t'\in T'$ and $ r\in R$. Since $T\in \Ore_l(R/ \ga )$, $\bt \br = \br_1\overline{t'}$  for some elements $t\in T'$ and $r_1\in R$. Then $s(tr-r_1t')=0$ for some element $s\in S$, and so $st\cdot r = sr_1\cdot t'$ where $st\in T'$, i.e.\ $T'\in \Ore_l(R)$.

(ii) $T'\in \Den_l(R)$: If $rt'=0$ for some elements $ r\in R$ and $t'\in T'$ then $\br \overline{t'}=0$, and so $\overline{t_1'}\br =0$ for some element $t_1'\in T'$. Now, $t_1'r\in \ga$. There exists an element $ s\in S$ such that $st_1'\cdot r=0$. Notice that $st_1'\in T'$ (since $s\in S\subseteq T'$). Therefore, $T'\in \Den_l(R)$, by (i).

(iii) $\ass (T')  =\gb$: If $t'r=0$ for some elements $ t'\in T'$ and $r\in R$ then $ \overline{t'}\br =0$, and so $\br \in \gb / \ga$ and $r\in \gb$. Then $\ass (T') \subseteq \gb$. Let $b\in \gb$. Then $ \overline{t'}\ob=0$ for some element $t'\in T'$ (since $\pi (T') = T\in \Den_l( R/ \ga , \gb / \ga )$). Now, $t'b\in \ga$. Hence, $st' \cdot b =0$ for some $s\in S$. Since $st'\in T'$, $b\in \ass (T')$. Then, $\ass (T') = \gb$. $\Box $

%$\noindent $

\begin{corollary}\label{b21Sep13}%\marginpar{b21Sep13}
Let $R$ be a ring and  $S\in \Den_l(R, \ga )$.  Then
\begin{enumerate}
\item $ S+\ga \in \Den_l(R, \ga )$.
\item $S+\ass (S) = S$ for all $S\in \maxDen_l(R)$.
\item $S+\ass (S) = S$ for all $S\in \maxDen (R)$.
\end{enumerate}
\end{corollary}

{\it Proof}. 1.    Let $\pi : R\ra R/ \ga$, $r\mapsto \br := r+\ga$. Since $T:= \pi (S)\in \Den_l(R/ \ga , 0)$ and $S+\ga = \pi^{-1}(T)\supseteq S$, we see that  $ S+\ga \in \Den_l(R, \ga )$, by Lemma \ref{a21Sep13}.

2 and 3. Statements 2 and 3 follow from statement 1 and the inclusion $S\subseteq S+\ass (S)$.   $\Box $

$\noindent $

{\bf The largest strong left quotient ring $Q_l^s(R)$ of a ring $R$}.

$\noindent $

{\it Definition}. The ring $Q_l^s(R):= T_l(R)^{-1}R$ is called the {\em largest strong left quotient ring} of the ring $R$.

$\noindent $

There are  exact sequences:

%\marginpar{SaQs1}
\begin{equation}\label{SaQs1}
0\ra  \glsR \ra R\ra Q_l^s(R), \;\; r\mapsto \frac{r}{1},
\end{equation}
%\marginpar{SaQs}
\begin{equation}\label{SaQs}
0\ra S_{l,0}(R)^{-1} \glsR \ra Q_l(R)
\stackrel{\th }{\ra} Q_l^s(R), \;\;\th (s^{-1}r)=s^{-1}r,
\end{equation}
where $s\in S_{l,0}(R)$ and $r\in R$ (if $\th (s^{-1} r) =0$ then $\frac{r}{1}=0$ in $Q_l^s(R)$, hence $r\in \glsR$, and so $\ker (\th  )= S_{l,0}^{-1}(R)^{-1} \glsR$).

For each $S\in \maxDen_l(R)$, there is a commutative diagram of ring homomorphisms:

%\marginpar{RQS}
\begin{equation}\label{RQS}
\xymatrix{
R\ar[rd]^{\s_S}\ar[r]^{\s } & Q^s_l(R)\ar[d]^{\s^s_S}\\
 & S^{-1} R}
\end{equation}
where, for  $r\in R$ and $s\in T_l(R)$,    $\s (r) = \frac{r}{1}$, $\s_S(r) = \frac{r}{1}$ and $\s_S^s(s^{-1}r) = s^{-1}r$. Clearly,
%\marginpar{krss}
\begin{equation}\label{krss}
\ker (\s^s_S) = T_l(R)^{-1}\ass (S).
\end{equation}
In more detail, if $\s_S^s(s^{-1}r) =0$ then $0= \frac{r}{1}\in S^{-1}R$, hence $r\in \ass (S)$, and so $\ker (\s_S^s ) = T_l(R)^{-1} \ass (S)$.
 Moreover, there is a canonical exact
sequence %\marginpar{RQS1}
\begin{equation}\label{RQS1}
0\ra T_l(R)^{-1} \gll_R \ra Q^s_l(R) \stackrel{\s^s }{\ra} \prod_{S\in \maxDen_l(R)}S^{-1}R, \;\; {\rm where}\;\; \s^s := \prod_{S\in \maxDen_l(R)}\, \s^s_S.
\end{equation}
In more detail, if $ \s^s (s^{-1}r)=0$ where $s\in T_l(R)$ and $r\in R$ then $S^{-1}R\ni \s_S^s (s^{-1}r) = s^{-1}r=0$ for all $S\in \maxDen_l(R)$, and so $r\in \ass (S)$ for all $S\in \maxDen_l(R)$, i.e.\ $r\in \cap_{S\in \maxDen_l(R)}\ass (S)=\gll_R$. Therefore, $\ker (\s^s) = T_l(R)^{-1}\gll_R$.

{\bf Criterion for $S_{l,0}(R) = T_l(R)$}. Recall that $S_{l,0}(R)\subseteq T_l(R)$ (Lemma \ref{c19Sep13}.(3)). The next lemma is a criterion for a ring $R$ to have the property that  $S_{l,0}(R)= T_l(R)$ (or/and that the  homomorphism  $\th : Q_l(R) \ra Q_l^s(R)$ is an isomorphism).
\begin{lemma}\label{a15Sep13}%\marginpar{a15Sep13}
Let $R$ be a ring. The following statements are equivalent.

\begin{enumerate}
\item $S_{l,0}(R) = T_l(R)$.
\item $\gll^s_R=0$.
\item $\th$ is an isomorphism.
\item $\th $ is a monomorphism.
\item $\th $ is a epimorphism and $S_{l,0}(R)+\gll^s_R\subseteq S_{l,0}(R)$.
\end{enumerate}
\end{lemma}

{\it Proof}. The following implications are obvious (see the exact sequence (\ref{SaQs})): $1\Leftrightarrow 2 \Leftrightarrow 4$, $ 1\Rightarrow 3\Rightarrow 4$ and $3\Rightarrow 5$.

$(5\Rightarrow 2 ) $ It suffices  to show that $T_l:= T_l(R)\subseteq S_l := S_{l,0}(R)$. Let $t\in T_l$. Then $Q_l^s(R) \ni t^{-1} = \th (s^{-1}r)$ for some elements $ s\in  S_l$ and $r\in R$. Hence,  $rt = s+a=:s'$ for some element $a\in \gll^s_R$.  Since, by the assumption,  $S_l+\gll^s_R\subseteq S_l$, we have $s'\in S_l$. Hence, $\ker (t\cdot ) =0$ (where $t\cdot : R\ra R$, $x\mapsto tx$) since $\ker (t\cdot ) \subseteq \ker (rt\cdot ) = \ker (s'\cdot ) =0$ as $s'\in S_l$. Thus $\gll^s_R=0$.  $\Box $

$\noindent $

{\bf Two more characterizations of the set $T_l(R)$}.
For  a ring $R$, let $\pCCR :=\{ c\in R\, | \, \ker (\cdot c) =0\}$ be the {\em set of left regular elements} of $R$ where $\cdot c : R\ra R$, $r\mapsto rc$.

$\noindent $

{\it Definition}. The sets
\begin{eqnarray*}
\CC^w_R&:=& \{ c\in R\, | \, \frac{c}{1}\in \CC_{S^{-1}R}\;\; {\rm for \; all}\;\; S\in \maxDen_l(R)\},   \\
 \pCCwR &:=& \{ c\in R\, | \, \frac{c}{1}\in {}'\CC_{S^{-1}R}\;\; {\rm for \; all}\;\; S\in \maxDen_l(R)\},
\end{eqnarray*}
are called the {\em sets of weak regular}  and {\em weak left regular elements} of $R$.

$\noindent $

The sets $\CC^w_R$ and $\pCCwR$ are multiplicative sets such that   $R^* \subseteq \CC^w_R\subseteq \pCCwR$.

\begin{lemma}\label{a22Sep13}%\marginpar{a22Sep13}
Let $R$ be a ring. Then
\begin{enumerate}
\item $\CC_R\subseteq \pCCR\subseteq \pCCwR$.
\item $\CLsl (R) \subseteq \CC^w_R$.
\end{enumerate}
\end{lemma}

{\it Proof}. 1. We have to show that $\pCCR \subseteq \pCCwR$. Given $c\in \pCCR$. Suppose that $c\not\in \pCCwR$, we seek a contradiction. Then there exist $S\in \maxDen_l(R)$ and $r\in R$ such that $\frac{r}{1}\frac{c}{1}=0$ and $\frac{r}{1}\neq 0$. Then $rc \in \ass (S)$, and so $src =0$ for some $s\in S$. Now, $sr=0$ (since $c\in \pCCR$). Therefore, $\frac{r}{1}=0$, a contradiction.

2. Statement 2 follows from Proposition \ref{a19Sep13}.(1). $\Box $

%$\noindent $

By Lemma \ref{a22Sep13}.(2),
%\marginpar{TLCp}
\begin{equation}\label{TLCp}
T_l(R)\subseteq \CLsl (R)  \subseteq \CC^w_R \subseteq \pCCwR .
\end{equation}
Two more characterizations of the set $T_l(R)$ are given below.
\begin{proposition}\label{b22Sep13}%\marginpar{b22Sep13}
\begin{enumerate}
\item The set $T_l(R)$ is the largest left denominator set  in the set $\pCCwR$.
\item The set $T_l(R)$ is the largest left denominator set  in the set $\CC^w_R$.
\end{enumerate}
\end{proposition}

{\it Proof}. 1.  Given $T\in \Den_l(R)$ such that $T\subseteq \pCCwR$. We have to show that $T\subseteq T_l(R)$.

(i) {\em For all} $S\in \maxDen_l(R)$, $ST\in \Den_l(R)$: By Lemma \ref{a14Sep13}.(2),  we have to show that  $0\not\in ST$. Suppose that $0\in ST$ for some $S\in \maxDen_l(R)$, we seek a contradiction. Then $s_1t_1\cdots s_nt_n =0$ for some elements $s_i\in S$ and $t_i\in T$, and in the ring $S^{-1}R$, $s_1t_1\cdots s_nt_n/1= 0$. Now,  $s_1t_1\cdots s_nt_n/1=0$ implies
 $s_1t_1\cdots s_n/1=0$ (since $t_n\in \pCCwR$) implies $s_1t_1\cdots s_{n-1}t_{n-1}/1=0$ (since $s_n$ is a unit in $S^{-1}R$). Continue in this way we obtain that $ s_1/1=0$, a contradiction.

(ii) $T\subseteq\CLsl (R)$. By (i), for all $S\in \maxDen_l(R)$, $S\subseteq ST$, hence $S= ST$ (by the maximality of $S$), and so $T\subseteq ST= S$. Therefore, $ T\subseteq \cap_{S\in \maxDen_l(R)} S= \CLsl (R)$.

(iii) $T\subseteq T_l(R)$, by Proposition \ref{b19Sep13}.

2. By (\ref{TLCp}), $\pCCwR\supseteq \CLsl (R) \supseteq T_l(R)$. Now, statement 2 follows from statement 1.  $\Box $

$\noindent $

%%%%%%%%%%%%%%%%%% SECTION 3  %%%%%%%%%%%%%%%%%%%%%%%%

\section{The largest strong left quotient ring of a ring and its properties}\label{TSLQRP}%\marginpar{TSLQRP}

The aim of this section is to  prove Theorem \ref{21Sep13} and Theorem \ref{B21Sep13}, and to  give a criterion for a ring $R$ to have  a semisimple  strong left quotient ring (Theorem \ref{Y20Sep13}).

Let us collect/prove the results (Lemma \ref{a6Jul10}, Lemma \ref{c21Sep13} and Proposition \ref{A8Dec12}) that are used in the proofs of these theorems. Let $S,T\in \Den_l(R)$. The denominator set $T$ is called $S$-{\em
saturated} if $sr\in T$, for some $s\in S$ and $r\in R$, then
$r\in T$, and if $r's'\in T$, for some $s'\in S$ and $r'\in R$,
then  $r'\in T$.

\begin{lemma}\label{a6Jul10}\cite{larglquot}%\marginpar{a6Jul10}
%\cite{larglquot}
$\,$ Let $S\in \Den_l(R, \ga )$, $\pi : R\ra R/\ga$, $a\mapsto a+\ga$,
and  $\s : R\ra S^{-1}R$, $ r\mapsto r/1$.
\begin{enumerate}
\item  Let $T\in \Den_l(S^{-1}R, 0)$ be such that $\pi (S), \pi
(S)^{-1}\subseteq T$. Then $T':= \s^{-1} (T) \in \Den_l(R, \ga )$,
$T'$ is $S$-saturated, $T=\{ s^{-1}t'\, | \, s\in S, t'\in T'\}$,
and $S^{-1} R\subseteq T'^{-1}R= T^{-1}R$. \item $\pi^{-1} (S_0(R/
\ga )) = S_\ga (R)$, $ \pi (S_\ga (R)) = S_0(R/\ga ))$ and $ Q_\ga
(R)=S_\ga (R)^{-1}R = Q_l(R/\ga )$.
\end{enumerate}
\end{lemma}

The next lemma shows that there is a canonical bijection between the sets of maximal left denominator sets of the rings $R$ and $R/ \glsR$.
\begin{lemma}\label{c21Sep13}%\marginpar{c21Sep13}
Let $R$ be a ring,
 % $\ga$ be an ideal of $R$ such that $\ga \subseteq \gll_R$ (for example, $\ga =\gll^s_R, \gll_R$) and
    $ \pi : R\ra R/ \glsR$, $r\mapsto \br := r+\glsR$. Then the map
$$ \maxDen_l (R) \ra \maxDen_l(R/\glsR ) , \;\; S\mapsto \pi (S),$$ is a bijection with the inverse $ T\mapsto \pi^{-1}(T)$.
\end{lemma}

{\it Proof}. (i) {\rm For all $S\in \maxDen_l(R)$, $\pi (S)\in \maxDen_l(R/ \glsR )$}:  Since $\glsR \subseteq \gll_R\subseteq \ass (S)$, $\pi (S)\in \Den_l(R/ \glsR , \ass (S)/ \glsR )$. Suppose that $T\in \Den_l(R/ \glsR )$ with $\pi (S)\subseteq T$, we have to show that $\pi (S)=T$. Let $\pi': R/ \glsR \ra R/\ass (S)$, $r+\glsR\mapsto r+\ass (S)$. Since $\ass (\pi (S)) \subseteq \ass (T)$,
$$ \pi'(T) \in \Den_l(R/ \ass (S), \ass (T) / \ass (\pi (S))).$$ Using Lemma \ref{a21Sep13} in the situation when $S\in \Den_l(R, \ass (S))$, $\pi'':R\ra R/ \ass (S)$, $r\mapsto r+\ass (S)$, and $\pi'(T)\in \Den_l( R/ \ass (S), \ass (T) / \ass (\pi (S))\simeq \gb / \ass (S))$ where $\gb = \pi^{-1} (\ass (T))$, we conclude that $$ S\subseteq \pi^{-1} \pi (S) \subseteq \pi^{-1}(T)= \pi''^{-1} (\pi'(T))\in \Den_l(R).$$ Therefore, $S= \pi^{-1} (T)$, by the maximality of $S$. Hence, $\pi (S) = \pi \pi^{-1} (T)=T$, as required.

(ii) {\em For all $T\in \maxDen_l(R/\glsR )$, $\pi^{-1}(T)\in \maxDen_l(R)$ }: Since $T_l:=T_l(R)\in \Den_l(R, \glsR )$, we have $\pi (T_l)\in \Den_l(R/ \glsR , 0)$. We claim that
$$ 0\not\in T\pi (T_l)$$
where $T\pi (T_l)$ is the submonoid of $(R/ \glsR , \cdot )$ generated by $T$ and $\pi (T_l)$.
Suppose that $0\in T\pi (T_l)$, i.e.\ $t_1s_1\cdots t_ns_n=0$ for some elements $ t_i\in T$ and $s_i\in \pi (T_l)$, we seek a contradiction. Then  $t_1s_1\cdots t_n=0$ (since $s_n\in \pi (T_l)$ and $\ass (\pi (T_l))=0$) and so
$$t_n' t_1s_1\cdot t_2s_2\cdots t_{n-1}s_{n-1}=0$$ for some element $t_n'\in T$  (since $t_n\in T$ and $T\in \Den_l(R/ \glsR )$). Repeating the same argument $n-1$ more times we obtain elements $ t_2', \ldots , t_{n-1}'\in T$ such that $T\ni t_2't_3'\cdots t_n't_1=0$, a contradiction.

Since $T\in \maxDen_l(R/ \glsR )$ and $0\not\in T\pi (T_l)$, we must have $T\subseteq T\pi (T_l) \in \Den_l(R/ \glsR )$, by Lemma \ref{a14Sep13}.(2). Therefore, $T= T\pi (T_l)$ (by the maximality of $T$)  and so
$$\pi (T_l)\subseteq T.$$  Using Lemma \ref{a21Sep13} in the situation when $S:= T_l(R)\in \Den_l(R, \ga := \glsR )$, $\pi : R\ra R/\glsR $ and $T\in \Den_l( R/\glsR , \ass (T))$, we conclude  that $ T_l(R)\subseteq \pi^{-1}(T)\in \Den_l(R)$. To finish the proof of (ii) we have to show that if $\pi^{-1}(T)\subseteq S'$ for some $S'\in \maxDen_l(R)$ than $\pi^{-1}(T)= S'$. The inclusion $\pi^{-1}(T)\subseteq S'$ implies the inclusion $T=\pi\pi^{-1} (T)\subseteq \pi (S')$. By (i) and $T\in \maxDen_l(R/ \glsR )$, $T= \pi (S')$. Therefore,
$$ \pi^{-1} (T) = \pi^{-1}\pi (S') = S'+\glsR = S'$$ since $\glsR\subseteq \gll_R\subseteq \ass (S')$ and $S'+\ass (S')= S'$, by Corollary \ref{b21Sep13}.(2).

(iii) {\em For all $S\in \maxDen_l(R)$, $\pi^{-1}\pi (S) = S$}: $ \pi^{-1}\pi (S) = S+\ass (S)= S$, by Corollary \ref{b21Sep13}.(2).

(iv) {\em For all $T\in \maxDen_l(R/\glsR )$, $\pi \pi^{-1} (T) = T$}: Trivial.

The proof of the lemma is complete. $\Box $

$\noindent $

{\bf A bijection between $\maxDen_l(R)$ and $\maxDen_l(Q_l(R))$}.
\begin{proposition}\label{A8Dec12}{\cite[Proposition 2.10]{Bav-Crit-S-Simp-lQuot}}%\marginpar{A8Dec12}
$\,$ Let $R$ be a ring, $S_l$ be the  largest regular left Ore set of the ring $R$, $Q_l:= S_l^{-1}R$ be the largest left quotient ring of the ring $R$, and $\CC$ be the set of regular elements of the ring $R$. Then
\begin{enumerate}
\item $S_l\subseteq S$ for all $S\in \maxDen_l(R)$. In particular,
$\CC\subseteq S$ for all $S\in  \maxDen_l(R)$ provided $\CC$ is a
left Ore set. \item Either $\maxDen_l(R) = \{ \CC \}$ or,
otherwise, $\CC\not\in\maxDen_l(R)$. \item The map $$
\maxDen_l(R)\ra \maxDen_l(Q_l), \;\; S\mapsto SQ_l^*=\{ c^{-1}s\,
| \, c\in S_l, s\in S\},
$$ is a bijection with the inverse $\CT \mapsto \s^{-1} (\CT )$
where $\s : R\ra Q_l$, $r\mapsto \frac{r}{1}$, and $SQ_l^*$ is the
sub-semigroup of $(Q_l, \cdot )$ generated by the set  $S$ and the
group $Q_l^*$ of units of the ring $Q_l$, and $S^{-1}R= (SQ_l^*)^{-1}Q_l$.
    \item  If $\CC$ is a left Ore set then the map (where $Q=Q_{l,cl}(R)$) $$ \maxDen_l(R)\ra \maxDen_l(Q), \;\; S\mapsto SQ^*=\{ c^{-1}s\,
| \, c\in \CC, s\in S\}, $$ is a bijection with the inverse $\CT
\mapsto \s^{-1} (\CT )$ where $\s : R\ra Q$, $r\mapsto
\frac{r}{1}$, and $SQ^*$ is the sub-semigroup of $(Q, \cdot )$
generated by the set  $S$ and the group $Q^*$ of units of the ring
$Q$, and $S^{-1}R= (SQ^*)^{-1}Q$.
\end{enumerate}
\end{proposition}

$\noindent $

{\bf Proof of Theorem \ref{21Sep13}}. 1. Let $T_l:= T_l(R)$ and $\CS := S_{l,\glsR }(R) $. Since $T_l\in \Den_l(R, \glsR )$, we have the inclusion $ T_l\subseteq \CS$. It remains to show that $T_l\supseteq \CS$. By Proposition  \ref{A8Dec12}.(1), and Lemma \ref{c21Sep13}, $S_{l,0}(R/\glsR ) \subseteq \pi (S)$ for all $S\in \maxDen_l(R)$. Therefore,
 \begin{eqnarray*}
\CS & =& \pi^{-1}(S_{l,0}(R/\glsR ))\;\;\; ({\rm Lemma}\;\; \ref{a6Jul10}.(2))\\
 &\subseteq & \pi^{-1} \left(\bigcap_{T\in \maxDen_l(R/\glsR )}T\right)\;\;\; ({\rm Proposition }\;\; \ref{A8Dec12}.(1))\\
  &= & \pi^{-1} \left(\bigcap_{S\in \maxDen_l(R )}\pi (S)\right)\;\;\; ({\rm Lemma }\;\; \ref{c21Sep13})\\
 &=& \bigcap_{S\in \maxDen_l(R )}\pi^{-1}\pi (S)=\bigcap_{S\in \maxDen_l(R )}(S+\glsR )\\
 &=&\bigcap_{S\in \maxDen_l(R )}S\;\;\; ({\rm by\;  Corollary }\;\; \ref{b21Sep13}.(2)\;\; {\rm and}\;\; \glsR\subseteq \ass (S))\\
 &=& \CL_l(R).
 \end{eqnarray*}
Hence, $ \CS \subseteq T_l$, by the maximality of $T_l$.

2. By statement 1, $Q_l^s(R)=T_l(R)^{-1}R=S_{l, \glsR}(R)^{-1}R= Q_{l , \glsR } (R)$. By Lemma \ref{a6Jul10}.(2), $Q_{l, \glsR } (R) \simeq Q_l(R/ \glsR )$.

3 and 4.  By Theorem \ref{4Jul10}.(1), $S_{l,0}(R/ \glsR )= R/ \glsR \cap Q_l(R/ \glsR )^*$. Now,
\begin{eqnarray*}
 T_l(R)&\stackrel{{\rm st.  1}}{=} &S_{l, \glsR }(R)=\pi^{-1} (S_{l,0}(R/ \glsR ))\;\;\;\; ({\rm Lemma }\; \ref{a6Jul10}.(2))\\
  &=& \pi^{-1} (R/ \glsR \cap Q_l(R/ \glsR )^*) \\
 &= & \pi^{-1} (R/ \glsR \cap Q_l^s(R/ \glsR )^*)\;\; ({\rm by \; statement\; 2}) \\
 &=&\s^{-1} (Q_l^s(R/ \glsR )^*).
\end{eqnarray*}
5. By Theorem \ref{4Jul10}.(3),
$$ Q_l(R/ \glsR )^* = \{ s^{-1}t\, | \, s,t\in S_{l,0}(R/ \glsR )\}. $$
By statements 1, 2 and 4, $Q_l^s(R)^* =  \{ s^{-1}t\, | \, s,t\in T_l(R )\}$. $\Box $

$\noindent $

{\bf  Proof of Theorem \ref{B21Sep13}}. 6. Statement 6 follows at once from Lemma \ref{c21Sep13}:
 \begin{eqnarray*}
 \pi^{-1}(\CLsl (R/ \glsR ))&=&\pi^{-1} \left( \bigcap_{T\in \maxDen_l(R/ \glsR )}T\right)=\bigcap_{T\in \maxDen_l(R/ \glsR )}\pi^{-1}(T)\\
 &=& \bigcap_{S\in \maxDen_l(R)}S =\CLsl (R).
\end{eqnarray*}
Hence, $\pi (\CLsl (R))= \pi \pi^{-1} (\CLsl (R/ \glsR ))= \CLsl (R/ \glsR ).$

2 and 4. Since $\pi (T_l(R))\in \Den_l(R/ \glsR , 0)$ and $\pi (T_l(R))\subseteq \pi ( \CLsl (R)) = \CLsl (R/ \glsR )$, by statement 6, we see that $\pi (T_l(R)) \subseteq T_l(R/ \glsR )$, by the maximality of $T_l(R/ \glsR )$.

Conversely, by Lemma \ref{c21Sep13}, $T_l(R/ \glsR ) \subseteq \pi (S)$ for all $S\in \maxDen_l(R)$. Therefore,
$$\pi^{-1} (T_l(R/ \glsR ))\subseteq \pi^{-1}\pi (S) = S\;\; {\rm  for \; all }\;\;  S\in \maxDen_l(R)\;\; ({\rm Lemma}\; \ref{c21Sep13}),$$
 hence $\pi^{-1} (T_l(R/ \glsR ))\subseteq \bigcap_{S\in \maxDen_l(R)}S=  \CLsl (R)$. Since $T_l(R/ \glsR )\in \Den_l(R/ \glsR )$, $T_l(R) \in \Den_l(R, \glsR )$ and $ \pi ( T_l(R))\subseteq T_l(R/ \glsR )$, we obtain that $T_l(R)\subseteq \pi^{-1} (T_l(R/ \glsR ))\in \Den_l(R)$, by Lemma \ref{a21Sep13}. Hence,
 $$\pi^{-1} (T_l(R/ \glsR ))= T_l(R),$$ by the maximality of $T_l(R)$. Hence,
 $$T_l(R/ \glsR )=\pi\pi^{-1} (T_l(R/ \glsR ))= \pi (T_l(R))\in \Den_l(R/ \glsR , 0).$$ This means that statements 2 and  4 hold.

 3.   Statement 3 follows from statement 4 and  Theorem \ref{21Sep13}.(1).
%\begin{eqnarray*}
%\pi^{-1} (T_l(R/ \glsR ))&=& \pi^{-1} (S_{l,0}(R/ \glsR )) \\
% &=& S_{l, \glsR}(R) \;\;\; {\rm by \; Lemma}\; \ref{a6Jul10}.(2)\\
% &=& T_l(R) \;\;\;\;\; \;\; {\rm by \; Theorem }\; \ref{21Sep13}.(1).
%\end{eqnarray*}

5. Let $\overline{\s} : R/ \glsR \ra Q_l(R/ \glsR )$, $\br \mapsto \frac{\br}{1}$. By Proposition \ref{A8Dec12}.(3), the map
%\marginpar{mDRa}
\begin{equation}\label{mDRa}
\maxDen_l(R/ \glsR ) \ra \maxDen_l(Q_l(R/ \glsR )), \;\; S\mapsto SQ_l(R/ \glsR )^*= \{ c^{-1} s\, | \, c\in S_{l, 0} (R/ \glsR ) , \; s\in S\},
\end{equation}
is a bijection with the inverse $T\mapsto \overline{\s}^{-1}(T) = T\cap R/ \glsR $. Let $\CT := T_l(Q_l(R/ \glsR ))$. Then $Q_l(R/ \glsR )^* \subseteq \CT$, by Lemma \ref{c19Sep13}.(3).
%In more detail,
%\begin{eqnarray*}
%Q_l(R/ \glsR )^* & \in & \Den_l(Q_l(R/ \glsR ),0) \\
% Q_l(R/ \glsR )^* & \subseteq  & \CLsl (Q_l(R/ \glsR )\;\;\;\;  {\rm and}\\
% Q_l(R/ \glsR )^*\CF  & \in  & \Den_l(Q_l(R/ \glsR )\;\;\;\;  {\rm and} \\
% Q_l(R/ \glsR )^*\CF  & \subseteq  & \CLsl (Q_l(R/ \glsR )
%\end{eqnarray*}
%hence $\CT \subseteq Q_l(R/ \glsR )^* \CT \subseteq \CT$, by the maximality of $\CT$.
 By (\ref{mDRa}), $\CT \subseteq SQ_l(R/ \glsR )^*$ for all $S\in \maxDen (R/ \glsR )$. Hence,
\begin{eqnarray*}
\CT \cap R/ \glsR  &\subseteq  & \bigcap_{S\in \maxDen_l(R/  \glsR )} R/ \glsR \cap SQ_l(R/ \glsR )^*\\
 &=&   \bigcap_{S\in \maxDen_l(R/ \glsR )}S\;\;\;\; ({\rm since}\;\; S =R/ \glsR \cap SQ_l(R/ \glsR )^*, \;\; {\rm by\; } \; (\ref{mDRa}))\\
 &=& \CLsl (R/ \glsR ).
\end{eqnarray*}
Also,
\begin{eqnarray*}
\CT \cap  R/ \glsR  &\supseteq  & R/ \glsR \cap Q_l(R/ \glsR )^*\;\; ({\rm since}\;\; \CT \supseteq Q_l(R/ \glsR)^*)\\
  &=& S_{l,0}(R/ \glsR ) \;\;\; {\rm (Theorem \; \ref{4Jul10}.(1))}\\
 &=& T_l(R/ \glsR ), \;\;\; {\rm (by \; statement\;\; 3)}.
\end{eqnarray*}
Applying Lemma \ref{a6Jul10} to the case where $S= S_{l,0}(R/ \glsR )\in \Den_l(R/ \glsR , \ga :=0 )$ and $T= \CT$ (notice that $Q_l(R/ \glsR )^*\subseteq \CT$) we see that
$$\CT \cap R/ \glsR  = \overline{\s}^{-1} (\CT ) \in \Den_l(R/ \glsR ).$$ This fact together with the inclusions (see above) $T_l(R/ \glsR ) \subseteq \CT \cap  R/ \glsR  \subseteq \CLsl (R/ \glsR )$ implies the equality
%\marginpar{TlRaR}
\begin{equation}\label{TlRaR}
T_l(R/ \glsR )=\CT \cap R/ \glsR ,
\end{equation}
by the maximality of $T_l(R/ \glsR )$. The inclusion $Q_l(R/ \glsR )^*\subseteq \CT$ and the fact that $ Q_l(R/ \glsR )^* =\{ s^{-1}t\, | \, s,t\in S_{l,0}(R/ \glsR )\}$ (Theorem \ref{4Jul10}.(3)) imply that
$$ \CT = \{ s^{-1} r\, | \, s\in S_{l,0}(R/ \glsR ); \; r\in \CT \cap  R/ \glsR  \}.$$
In more detail, if $s^{-1}r\in \CT$  for some $s\in S_{l,0}(R/ \glsR )\subseteq \CT$ and $\frac{r}{1}\in R/ \glsR$ then $\frac{r}{1}\in s\CT \subseteq \CT \CT = \CT$, and so $r\in \CT \cap R/ \glsR$.

By (\ref{TlRaR}) and $T_l(R/ \glsR ) = S_{l,0}(R/ \glsR )$ (statement 3),
$$ \CT = \{ s^{-1} r\, | \, s\in T_l(R/ \glsR ), \;\; r\in T_l(R/ \glsR )\} =Q_l(R/ \glsR )^*.$$Hence, $\ass (\CT )=\ass (Q_l(R/ \glsR )^*) =0$.

1. Statement 1 follows from statement 5.

7. By statement 4 and Theorem \ref{21Sep13}.(2), $Q_l^s(R/ \glsR ) \simeq Q_l((R/ \glsR )/ \gll^s_{R/ \glsR })=  Q_l(R/ \glsR )$. $\Box $

$\noindent $

{\bf Necessary and sufficient conditions for  $Q_l(R)$ to be  a
semi-simple ring}. A ring $Q$ is called a {\em ring of quotients}
if every element $c\in \CC_Q$ is invertible. A subring $R$ of a
ring of quotients $Q$ is called a {\em left order} in $Q$ if
$\CC_R$ is a left Ore set and $\CC_R^{-1}R=Q$. A ring $R$ has {\em
finite left rank} (i.e. {\em finite left uniform dimension}) if
there are no infinite direct sums of nonzero left ideals in $R$.

The next theorem gives an answer to the question of when $Q_l(R)$
is a semi-simple ring. Theorem \ref{5Jul10} is the key result in the proof of Theorem \ref{Y20Sep13}.

\begin{theorem}\label{5Jul10}\cite{larglquot}%\marginpar{5Jul10}
$\,$ The following properties of a ring $R$ are equivalent.
\begin{enumerate}
\item  $Q_l(R)$ is a semi-simple ring. \item $Q_{cl}(R)$  exists
and is a semi-simple ring. \item $R$ is a left order in a
semi-simple ring. \item $R$ has finite left rank, satisfies the
ascending chain condition on left annihilators and is a semi-prime
ring. \item A left ideal of $R$ is essential iff it contains a
regular element.
\end{enumerate}
If one of the equivalent conditions hold then $S_0(R) = \CC_R$ and
$Q_l(R) = Q_{cl}(R)$.
\end{theorem}

%$\noindent $

{\bf  Proof of Theorem \ref{Y20Sep13}}. $(2\Leftrightarrow 4)$  This is the Goldie's Theorem for the ring $R/ \glsR$.

 $(1\Leftrightarrow 3\Leftrightarrow 4)$  By Theorem \ref{21Sep13}.(2), $Q_l^s(R)\simeq Q_l(R/ \glsR )$. Theorem \ref{5Jul10} implies that  $(1\Leftrightarrow 3\Leftrightarrow 4)$ and that $ Q_l(R/ \glsR ) \simeq Q_{l,cl}(R/ \glsR )$ and $S_{l,0}(R/ \glsR ) = \CC_{R/ \glsR}$. By Theorem \ref{B21Sep13}.(3), $S_{l,0}(R/ \glsR )= T_l(R/ \glsR )$. Then, by Theorem \ref{B21Sep13}.(2), $T_l(R)=\pi^{-1} (T_l(R/ \glsR ))=\pi^{-1} (S_{l,0}(R/ \glsR ))$. By Theorem \ref{B21Sep13}.(2),   $T_l(R/ \glsR ) = \pi (T_l)=\pi\pi^{-1} (S_{l,0}(R/ \glsR )) =  S_{l,0}(R/ \glsR )=\CC_{R/ \glsR}$. $\Box $

$\noindent $

{\bf The maximal left denominator sets  of a finite direct product of rings}.
\begin{theorem}\label{c26Dec12}{\cite[Theorem 2.9]{Bav-Crit-S-Simp-lQuot}}%\marginpar{c26Dec12}
$\,$ Let $R=\prod_{i=1}^n R_i$ be a direct product of rings $R_i$. Then
for each $i=1, \ldots , n$, the map
%\marginpar{1aab1}
\begin{equation}\label{1aab1}
\maxDen_l(R_i) \ra \maxDen_l(R), \;\; S_i\mapsto R_1\times\cdots \times S_i\times\cdots \times R_n,
\end{equation}
is an injection. Moreover, $\maxDen_l(R)=\coprod_{i=1}^n \maxDen_l(R_i)$ in the sense of (\ref{1aab1}), i.e.
$$ \maxDen_l(R)=\{ S_i\, | \, S_i\in \maxDen_l(R_i), \; i=1, \ldots , n\},$$
$S_i^{-1}R\simeq S_i^{-1}R_i$, $\ass_R(S_i)= R_1\times \cdots \times \ass_{R_i}(S_i)\times\cdots \times R_n$.
%The core of the left denominator set $S_i$ in $R$ coincides with the core $S_{i,c}$ of the left denominator set $S_i$ in $R_i$, i.e.\
%$$(R_1\times\cdots \times S_i\times\cdots \times R_n)_c=0\times\cdots \times S_{i,c}\times\cdots \times 0.$$
\end{theorem}

{\bf Proof of Theorem \ref{27Sep13}}. $(1\Leftrightarrow 2)$  These implications are Theorem \ref{5Jul10}.(1,2).

$(3\Rightarrow 1)$ The implication follows from Theorem \ref{Y20Sep13}. In particular, $Q_{l,cl}(R)\simeq Q_l^s(R)$.
 %By the statement (c), $T(R)\subseteq \CC_R$. Therefore, $T_l(R)=\CC_R$, by the statement (b). Thus, $Q_{l,cl}(R) = Q^s_l(R)$ is a semisimple  ring, by the statement (a).

 $(1\Rightarrow 3)$ (i) $\CLsl (R) = \CC_R$: Since $Q_{l,cl}(R)$ is a semisimple ring, $\CC_R= S_{l,0}(R)\subseteq \CLsl (R)$, by (\ref{SlinSL}). It remains to show that $\CLsl (R)\subseteq \CC_R$. The ring $Q_{l,cl}(R) = \prod_{i=1}^n R_i$ is a semisimple ring where $R_i$ are simple Artinian rings (i.e. matrix rings over division rings). Clearly, $\maxDen_l(R_i)=\{ R_i^*\}$ where $R_i^*$ is the group of units of $R_i$. By Theorem \ref{c26Dec12}, $\maxDen_l(R)=\{ \CS_1, \ldots , \CS_n\}$ where $\CS_i=R_1\times\cdots \times R_{i-1}\times R_i^*\times R_{i+1}\times \cdots \times R_n$ for $i=1, \ldots , n$. The map
 $$ \s : R\ra Q_{l,cl } (R) = \prod_{i=1}^n R_i, \;\; r\mapsto \frac{r}{1}=(r_1, \ldots , r_n),$$ is a monomorphism. Then an element $r\in R$ is regular iff the element $\frac{r}{1}=(r_1, \ldots , r_n) \in Q_{l, cl}(R)=\prod_{i=1}^n R_i$ is regular  iff $r_1\in \CC_{R_i}=R_1^*, \ldots , r_n\in \CC_{R_n}= R_n^*$ iff $\frac{r}{1}\in \bigcap_{i=1}^n \CS_i$ iff
 $$r=\s^{-1} (r)\in \s^{-1} (\bigcap_{i=1}^n \CS_i)=\bigcap_{i=1}^n \s^{-1}(\CS_i) = \bigcap_{S\in \maxDen_l(R)} S= \CLsl (R),$$
by Proposition \ref{A8Dec12}.(3) and Theorem \ref{c26Dec12}.

 (ii) $T_l(R) = \CC_R$: By (\ref{SlinSL}), $\CC_R=S_{l,0}(R) \subseteq T_l(R) \subseteq \CLsl (R) = \CC_R$.

 (iii) $Q_{l,cl}(R) = Q_l^s(R)$ is a semisimple ring, by (ii) and $\glsR = \ass (\CC_R) =0$. $\Box$

The operations $\CLsl (\cdot ):R\mapsto \CLsl (R)$, $T_l(\cdot) : R\mapsto T_l(R)$ and $ Q_l^s(\cdot ): R\mapsto  Q_l^s(R)$ commute with finite direct products as the next theorem shows.

\begin{theorem}\label{A21Sep13}%\marginpar{A21Sep13}
Let $R=\prod_{i=1}^n R_i$ be a direct product of rings. Then
\begin{enumerate}
\item $\CLsl (R) = \prod_{i=1}^n \CLsl (R_i)$
\item $T_l(R) = \prod_{i=1}^n T_l(R_i)$ and $ Q_l^s(R) \simeq \prod_{i=1}^nQ_l^s(R_i)$.
\end{enumerate}
\end{theorem}

{\it Proof}.  1. Statement 1 follows from Theorem \ref{c26Dec12}:
\begin{eqnarray*}
 \CLsl (R)&=& \bigcap_{S\in \maxDen_l(R)}S=\prod_{i=1}^n \bigcap_{S_i\in \maxDen_l(R_i)}S_i\;\; ({\rm Theorem }\;\; \ref{c26Dec12})\\
 &=&\prod_{i=1}^n\CLsl (R_i). \;\;\;
\end{eqnarray*}

2.  Let $T_l:= T_l(R)$ and $T = \prod_{i=1}^n T_l(R_i)$. We have to show that $T_l=T$. Clearly, $T\in \Den_l(R)$ and $ T\subseteq \prod_{i=1}^n \CLsl (R_i)=\CLsl (R)$, by statement 1. Therefore, $T\subseteq T_l$, by the maximality of $T_l$. It remains to show that $T_l\subseteq T$. Since $T_l = \prod_{i=1}^n T_i'$ for some $T_i'\in \Den_l(R_i)$ such that $T_i'\subseteq \CLsl (R_i)$ (by statement 1). Therefore, $T_i'\subseteq T_l(R_i)$ for $i=1, \ldots , n$, by the maximality of $T_l(R_i)$, and so $T_l\subseteq T$. $\Box$

\begin{lemma}\label{a27Nov12}%\marginpar{a27Nov12}
\cite{larglquot} Let $S\in \Den_l(R, \ga )$ and $T\in \Den_l(R, \gb )$ such that $ \ga \subseteq \gb$. Then
\begin{enumerate}
\item ${\rm r.ass} (ST)\subseteq \gb$ where ${\rm r.ass} (ST):=\{ r\in R\, | \, rc=0$ for some $c\in ST\}$.
\item $ST \in \Den_l(R, \gc )$ and $\gb \subseteq \gc$.
\end{enumerate}
\end{lemma}

{\bf Proof of Theorem \ref{A29Sep13}}. (i) $S_{l,\ga}(R)\subseteq S$ {\em for all} $S\in \maxDen_l(R)$: By Lemma \ref{a27Nov12}, $  S_{l,\ga}(R)S\in \Den_l(R)$ since $\ga \subseteq \glsR \subseteq \ass (S)$. Therefore, $S_{l,\ga}(R)\subseteq S_{l,\ga}(R)S = S$, by the maximality of $S$.

(ii) $S_{l,\ga}(R)\subseteq T_l(R)$: By (i), $S_{l,\ga}(R)\subseteq \CLsl (R)$, hence $S_{l,\ga}(R)\subseteq T_l(R)$, by the maximality of $T_l(R)$, Proposition \ref{b19Sep13}. Hence, there is a ring $R$-homomorphism $Q_{l,\ga}(R)\ra Q_l^s(R)$, $s^{-1}r\mapsto s^{-1}r$.  $\Box$

{\bf The largest strong  quotient ring of $Q_l(R)$}\label{LSQQlR}%\marginpar{LSQQlR}

For an arbitrary ring $R$,   Theorem \ref{E29Aug14}  establishes natural and tight connections between triples $T_l(R)$, $\gll_R^s$, $Q_l^s(R)$ and  $T_l(Q_l(R))$, $\gll_{Q_l(R)}^s$, $Q_l^s(Q_l(R))$. The applications of this theorem are given in Section \ref{EXAMPQ} where it is used in giving explicit descriptions of the triple  $T_l(R)$, $\gll_R^s$, $Q_l^s(R)$ for every ring $R$ such that its classical  left quotient ring $Q_{l,cl}(R)$ is a left Artinian ring, see Theorem \ref{F29Aug14} (notice that in this case $Q_{l,cl}(R)=Q_l(R)$ ({\cite[Corollary 2.10]{larglquot}}).

$\noindent $

{\bf Proof of Theorem \ref{E29Aug14}}. 1 and 2. Let $S_l=S_l(R)$ and $Q=Q_l(R)$. Since $S_l(R)\subseteq T_l(R)$ (Lemma \ref{c19Sep13}.(3)), the multiplicative submonoid $Q^*T_l(R)$ of $Q$ generated  by $Q^*$ and $T_l(R)$ belongs to $\Den_l(Q, S_l(R)^{-1}\gll^s_R)$. In view of the bijection between the sets $\maxDen_l(R)$ and $\maxDen_l(Q)$ given by Proposition \ref{A8Dec12}.(3), we have the inclusion $Q^*T_l(R)\subseteq \CL_l^s(Q)$ which immediately implies the inclusion $Q^*T_l(R)\subseteq T_l(Q)$. We identify the ring $R$ with its image in $Q$ via $\s : R\ra Q$, $r\mapsto \frac{r}{1}$.  By Lemma \ref{a6Jul10}.(1),  $T_l(Q)\cap R\in \Den_l(R, R\cap \gll^s_Q)$ and $T_l(Q)
=\{ s^{-1}t\, | \, s\in S_l(R), t\in T_l(Q)\cap R\}$. By Proposition \ref{A8Dec12}.(3), $T_l(Q)\cap R\subseteq \CL_l^s(R)$, and so $T_l(Q)\cap R \subseteq T_l(R)$. Now,
$$ T_l(R)\subseteq Q^*T_l(R)\cap R \subseteq T_l(Q)\cap R \subseteq T_l(R),$$
i.e. $T_l(R) = Q^*T_l(R)\cap R =T_l(Q)\cap R\in \Den_l(R, R\cap \gll^s_R)$. In particular, $\gll^s_R= R\cap \gll^s_Q$. By Lemma \ref{a6Jul10}.(1),
$$ Q^*T_l(R) = \{ s^{-1}t\, | \, s\in S_l(R), t\in T_l(R)=T_l(Q)\cap R\} = T_l(Q).$$ So, statement 1 is proven. Since $T_l(Q) = Q^*T_l(R)\in \Den_l(Q, S_l(R)^{-1}\ass (T_l(R)))$, we must have $\gll^s_Q= \ass (T_l(Q))= S_l(R)^{-1}\gll_R^s$. So, statement 2 is proven.

3. Let $Q=Q_l(R)$. Then
$$ Q_l^s(R)=T_l(R)^{-1}R\simeq (Q^*T_l(R))^{-1}S_l(R)^{-1}R\stackrel{{\rm st. \, 1}}{\simeq} T_l(Q)^{-1}Q\simeq Q_l^s(Q).\;\; \Box $$

%$\noindent $

%%%%%%%%%%%%%%%%%% SECTION 4  %%%%%%%%%%%%%%%%%%%%%%%%

\section{The largest strong  quotient ring of a ring}\label{T4SQR}%\marginpar{T4SQR}

In this section, the two-sided versions of the concepts appeared in Sections
 \ref{PRLM} and \ref{TSLQRP} are introduced: the {\em largest strong denominator set} $T(R)$, the {\em largest strong quotient ring} $Q^s(R)=T(R)^{-1}R$ and the {\em strong localization radical} $\glsRt $ (the subscript $t$ stands for `two-sided', i.e. `left and right'). All the results of the previous  sections  are true (with obvious adjustments)  for {\em left and right} Ore/denominator sets. For the analogous versions we state the corresponding results and the proofs are left for the reader as an exercise in the case when they are literally the same (with obvious modifications). The following notation is fixed:

\begin{itemize}

%\item   $\gn =\gn_R$ is the  prime radical of $R$ and $\nu\in \N \cup \{ \infty \}$ is its {\em nilpotency degree} ($\gn^\nu \neq 0$ but $\gn^{\nu +1}=0$);
%\item   $\bR := R/ \gn$ and $\pi: R\ra \bR$, $r\mapsto \br =r+\gn$;
%\item   $\OCC := \CC_{\bR}$ is the set of regular elements of the ring $\bR$ and $\bQ := \OCC^{-1}\bR$ is its left quotient ring;
%\item   $\CC':= \pi^{-1}(\OCC):=\{ c\in R\, \, | \, c+\gn \in \OCC\}$ and $Q':=\CC'^{-1}R$ (if it exists),
        \item $\Den (R, \ga )$ is the set of (left and right) denominator sets $S$ of $R$ with $\ass (S)=\ga$;
        \item $S_\ga=S_\ga (R)=S_{l,r, \ga }(R)$
 is the largest element of the poset $(\Den (R, \ga ),
\subseteq )$, i.e. the {\em largest denominator set in $R$ associated with} $\ga$,  and $Q_\ga (R):=Q_{l,r, \ga }(R):=S_{l,r,\ga}^{-1} R$ is  the
{\em largest (left and right) quotient ring associated with} $\ga$, $S_{l,r,\ga} $
exists,  \cite{larglquot};
\item $\maxDen (R)$ is the set of maximal  denominator sets of $R$ (it is a {\em non-empty} set, Lemma \ref{bb23Sep13}.(3));
    \item $\Ass (R):= \{ \ass (S)\, | \, S\in \Den (R)\}$.
\end{itemize}
The sets
$$ \CL (R):= \bigcup_{S\in \maxDen (R)}S\;\; {\rm and}\;\; \CN\CL (R):=R\backslash \CL (R)$$ are called the sets of {\em  localizable} and {\em  non-localizable elements} of $R$, respectively, and the intersection
$$\CL^s (R) :=\bigcap_{S\in \maxDen (R)} S$$
 is called the {\em set of strongly (or completely)   localizable elements} of $R$. Clearly, $\CL^s (R)$ is a multiplicative set and
 %\marginpar{TRsSl}
\begin{equation}\label{TRsSl}
 R^*\subseteq \CL^s(R)
\end{equation}
 since $R^*\subseteq S$ for all $S\in \maxDen (R)$, by Lemma \ref{a27Nov12}.(2). Similarly, by  Lemma  \ref{a27Nov12}.(2),
%\marginpar{TSlinSL}
\begin{equation}\label{TSlinSL}
S_{l,r,0}(R)\subseteq  \CL^s(R).
\end{equation}
%In particular, if the set $\CC_R$ of regular elements of the ring $R$ is a left Ore set then $\CC_R = S_{l,0}(R)$ and so
%%\marginpar{1SlinSL}
%\begin{equation}\label{1SlinSL}
%\CC_R \subseteq  \CL^s_l(R).
%\end{equation}
The next proposition is a characterization of the set $\CL^s(R)$ which says that the set $\CL^s (R)$ contains {\em precisely} the elements of the ring $R$ that are units in all maximal  localizations of $R$.

\begin{proposition}\label{Ta19Sep13}%\marginpar{Ta19Sep13}
Let $R$ be a ring. Then
\begin{enumerate}
\item $\CL^s(R) = \{ s\in R\,| \, \frac{s}{1}\in (S^{-1}R)^*$ for all $S\in \maxDen (R)\}$ where $(S^{-1}R)^*$ is the group of units of the ring $S^{-1}R^*$.
\item For all automorphisms $\s \in \Aut (R)$, $\s (\CL^s (R))= \CL^s (R)$.
\end{enumerate}
\end{proposition}

{\it Proof}. 1.  Let $\CR$ be the RHS of the equality. By the very definition of the set $\CL^s(R)$, $\CL^s(R) \subseteq \CR$. Conversely, let $s\in \CR$ and $\s_S :R\ra S^{-1}R$, $r\mapsto \frac{r}{1}$, where $S\in \maxDen (R)$. Then $s\in \s^{-1}_S((S^{-1}R)^*)= S$ for all $S\in \maxDen (R)$, \cite[Theorem 4.11.(2,3)]{larglquot}, hence $s\in \CL^s(R)$.

2. Obvious.   $\Box $

$\noindent $

{\bf The largest strong  denominator set $T(R)$}.
 The set $\Den^s (R):= \{ T\in \Den (R)\, | \, T\subseteq \CL^s_l(R)\}$ is a {\em non-empty} set since $ R^*\in \Den^s (R)$.  The elements of $\Den^s(R)$ are called the {\em strong  denominator sets} of $R$ and the rings $T^{-1}R$ where $T\in \Den^s (R)$ are called the {\em strong quotient rings} or the {\em  strong  localizations} of $R$. Proposition \ref{Tb19Sep13} shows that the set of maximal elements  ${\rm max.Den}^s(R)$ of the poset $(\Den^s(R), \subseteq )$  is a non-empty set. Moreover, it contains a single element. Namely,
 %\marginpar{TTlR}
\begin{equation}\label{TTlR}
T(R):=\bigcup_{S\in \Den^s (R)}S.
\end{equation}

\begin{proposition}\label{Tb19Sep13}%\marginpar{Tb19Sep13}
Let $R$ be a ring. Then ${\rm max.Den}^s(R)=\{T(R)\}$.
\end{proposition}

So, $T(R)$ is the {\em largest strong denominator set} of $R$.
 The ideal of the ring $R$,
$$ \glsRt := \ass (T(R))= \bigcup_{S\in \Den^s (R)}\ass (S)$$ is called the {\em strong  localization radical } of the ring $R$.

$\noindent $

  The ideal
$$ \gll_{R,t}:=\bigcap_{S\in \maxDen (R)} \ass (S)$$
is called the {\em (two-sided) localization radical} of $R$.    There is an exact sequence
$$ 0\ra \gll_{R,t}\ra R\stackrel{\s}{\ra} \prod_{S\in \maxDen (R)} S^{-1} R , \;\; \s =\prod_{S\in \maxDen (R)}\s_S, $$
where $\s_S:R\ra S^{-1}R$, $r\mapsto \frac{r}{1}$.
 The set $\Den^s (R)$ is invariant under the action of the group of automorphisms $\Aut (R)$ of the ring $R$.
\begin{lemma}\label{Tbb19Sep13}%\marginpar{Tbb19Sep13}
\begin{enumerate}
\item For all automorphisms $\s \in \Aut (R)$, $\s (T(R))= T(R)$, $\s (\gll_{R,t})= \gll_{R,t}$  and  $ \s (\glsRt)= \glsRt$.
\item $\glsRt\subseteq \gll_{R,t}$.
\end{enumerate}
\end{lemma}

%$\noindent $

\begin{lemma}\label{Tc19Sep13}%\marginpar{Tc19Sep13}
Let $R$ be a ring. Then
\begin{enumerate}
\item $T (R)+\glsRt\subseteq T(R)$.
\item $\CL^s(R)+\gll_{R,t}\subseteq \CL^s(R)$. In particular,  $T(R) + \gll_{R,t}\subseteq \CL^s(R)$.
\item $R^*\subseteq S_{l,r,0}(R)\subseteq T(R)$.
\end{enumerate}
\end{lemma}

{\bf The largest strong  quotient ring $Q^s(R)$ of a ring $R$}.
The ring $Q_l(R):= T_l(R)^{-1}R$ is called the {\em largest strong  quotient ring} of the ring $R$.  There are  exact sequences:

%\marginpar{TSaQs1}
\begin{equation}\label{TSaQs1}
0\ra  \glsRt \ra R\ra Q^s(R), \;\; r\mapsto \frac{r}{1},
\end{equation}
%\marginpar{TSaQs}
\begin{equation}\label{TSaQs}
0\ra S_{l,r,0}(R)^{-1} \glsRt \ra Q(R):=S_{l,r,0}(R)^{-1}R
\stackrel{\th }{\ra} Q^s(R), \;\;\th (s^{-1}r)=s^{-1}r,
\end{equation}
where $s\in S_{l,r,0}(R)$ and $r\in R$ (if $\th (s^{-1} r) =0$ then $\frac{r}{1}=0$ in $Q^s(R)$, hence $r\in \glsRt$, and so $\ker (\th  )= S_{l,r,0}^{-1}(R)^{-1} \glsRt$) and $Q(R)$ is the {\em largest (two-sided) quotient ring} of $R$, \cite{larglquot}.

For each $S\in \maxDen (R)$, there is a commutative diagram of ring homomorphisms:

%\marginpar{TRQS}
\begin{equation}\label{TRQS}
\xymatrix{
R\ar[rd]^{\s_S}\ar[r]^{\s } & Q^s(R)\ar[d]^{\s^s_S}\\
 & S^{-1} R}
\end{equation}
where, for  $r\in R$ and $s\in T(R)$,    $\s (r) = \frac{r}{1}$, $\s_S(r) = \frac{r}{1}$ and $\s_S^s(s^{-1}r) = s^{-1}r$. Clearly,
%\marginpar{Tkrss}
\begin{equation}\label{Tkrss}
\ker (\s^s_S) = T(R)^{-1}\ass (S).
\end{equation}
In more detail, if $\s_S^s(s^{-1}r) =0$ then $0= \frac{r}{1}\in S^{-1}R$, hence $r\in \ass (S)$, and so $\ker (\s_S^s ) = T(R)^{-1} \ass (S)$.
 Moreover, there is a canonical exact
sequence %\marginpar{TRQS1}
\begin{equation}\label{TRQS1}
0\ra T(R)^{-1} \gll_{R,t} \ra Q^s(R) \stackrel{\s^s }{\ra} \prod_{S\in \maxDen (R)}S^{-1}R, \;\; {\rm where}\;\; \s^s := \prod_{S\in \maxDen (R)}\, \s^s_S.
\end{equation}
In more detail, if $ \s^s (s^{-1}r)=0$ where $s\in T(R)$ and $r\in R$ then $S^{-1}R\ni \s_S^s (s^{-1}r) = s^{-1}r=0$ for all $S\in \maxDen (R)$, and so $r\in \ass (S)$ for all $S\in \maxDen (R)$, i.e.\ $r\in \cap_{S\in \maxDen (R)}\ass (S)=\gll_{R,t}$. Therefore, $\ker (\s^s) = T (R)^{-1}\gll_{R,t}$.

$\noindent $

{\bf Criterion for $S_{l,r,0}(R) = T(R)$}. Recall that $S_{l,r,0}(R)\subseteq T_l(R)$ (Lemma \ref{Tc19Sep13}.(3)). The next lemma is a criterion for a ring $R$ to have the property that  $S_{l,r,0}(R)= T(R)$ (or/and that the homomorphism $\th : Q_l(R) \ra Q^s(R)$ is an isomorphism).
\begin{lemma}\label{Ta15Sep13}%\marginpar{Ta15Sep13}
Let $R$ be a ring. The following statements are equivalent.

\begin{enumerate}
\item $S_{l,r,0}(R) = T(R)$.
\item $\gll^s_{R,t}=0$.
\item $\th$ is an isomorphism.
\item $\th $ is a monomorphism.
\item $\th $ is a epimorphism and $S_{l,r,0}(R)+\gll^s_{R,t}\subseteq S_{l,r,0}(R)$.
\end{enumerate}
\end{lemma}

%$\noindent $

{\bf Two more characterizations of the set $T(R)$}.
 The sets
\begin{eqnarray*}
\CC^w_{R,t}&:=& \{ c\in R\, | \, \frac{c}{1}\in \CC_{S^{-1}R}\;\; {\rm for \; all}\;\; S\in \maxDen (R)\},   \\
 \pCCwRt &:=& \{ c\in R\, | \, \frac{c}{1}\in {}'\CC_{S^{-1}R}\;\; {\rm for \; all}\;\; S\in \maxDen (R)\},
\end{eqnarray*}
are called the {\em sets of two-sided weak regular}  and {\em weak left regular elements} of $R$, respectively.

$\noindent $

The sets $\CC^w_{R,t}$ and $\pCCwRt$ are multiplicative sets such that   $R^* \subseteq \CC^w_{R,t}\subseteq \pCCwRt$.

\begin{lemma}\label{Ta22Sep13}%\marginpar{Ta22Sep13}
Let $R$ be a ring. Then
\begin{enumerate}
\item $\CC_R\subseteq \pCCR\subseteq \pCCwRt$.
\item $\CL^s (R) \subseteq \CC^w_{R,t}$.
\end{enumerate}
\end{lemma}

{\it Proof}. 1. We have to show that $\pCCR \subseteq \pCCwRt$. Given $c\in \pCCR$. Suppose that $c\not\in \pCCwRt$, we seek a contradiction. Then there exist $S\in \maxDen (R)$ and $r\in R$ such that $\frac{r}{1}\frac{c}{1}=0$ and $\frac{r}{1}\neq 0$. Then $rc \in \ass (S)$, and so $src =0$ for some $s\in S$. Now, $sr=0$ (since $c\in \pCCR$). Therefore, $\frac{r}{1}=0$, a contradiction.

2. Statement 2 follows from Proposition \ref{Ta19Sep13}.(1). $\Box $

%$\noindent $

By Lemma \ref{Ta22Sep13}.(2),
%\marginpar{TTLCp}
\begin{equation}\label{TTLCp}
T(R)\subseteq \CL^s(R)\subseteq \CC^w_{R,t}\subseteq \pCCwRt .
\end{equation}
The next proposition gives two more characterizations of the set $T(R)$.

\begin{proposition}\label{Tb22Sep13}%\marginpar{Tb22Sep13}
\begin{enumerate}
\item The set $T(R)$ is the largest left denominator set  in the set $\pCCwRt$.
\item The set $T(R)$ is the largest left denominator set  in the set $\CC^w_{R,t}$.
\end{enumerate}
\end{proposition}

The next lemma shows that there is a bijection between maximal denominator sets of the rings $R$ and $R/ \glsRt$.

\begin{lemma}\label{Tc21Sep13}%\marginpar{Tc21Sep13}
Let $R$ be a ring,
 % $\ga$ be an ideal of $R$ such that $\ga \subseteq \gll_R$ (for example, $\ga =\gll^s_R, \gll_R$) and
    $ \pi : R\ra R/ \glsRt$, $R\mapsto \br := r+\glsRt$. Then the map
$$ \maxDen (R) \ra \maxDen (R/\glsRt ) , \;\; S\mapsto \pi (S),$$ is a bijection with the inverse $ T\mapsto \pi^{-1}(T)$.
\end{lemma}

\begin{theorem}\label{T21Sep13}%\marginpar{T21Sep13}
Let $R$ be a ring, $\pi : R\ra R/ \glsRt$, $r\mapsto \br = r+\glsRt$; $\s : R\ra Q^s(R)$, $r\mapsto \frac{r}{1}$, and $Q^s(R)^*$ be the group of units of the ring $Q^s(R)$. Then
\begin{enumerate}
\item $T(R) = S_{l,r, \glsRt }(R)$.
\item $Q^s(R) = Q_{l,r,  \glsRt}(R) \simeq  Q_l(R/ \glsRt)$.
\item $T(R) = \s^{-1}(Q^s(R)^*)$.
\item $T(R) = \pi^{-1} (S_{l,r,0}(R/\glsRt))$.
\item $Q^s(R)^*= \{ s^{-1}t\, | \, s,t \in T(R)\}$.
\end{enumerate}
\end{theorem}

\begin{theorem}\label{TB21Sep13}%\marginpar{TB21Sep13}
We keep the notation of Theorem \ref{T21Sep13}. Then
\begin{enumerate}
\item  $Q^s(Q^s(R))=Q^s(R)$.
\item $T(R/ \glsRt) = \pi (T(R))$ and $ T(R) = \pi^{-1} (T ( R/ \glsRt))$.
\item $T(R/ \glsRt) =S_{l,r,0}(R/ \glsRt) $.
\item $\gll^s_{R/ \glsRt , t}=0$.
\item $T(Q^s(R)) = Q^s(R)^*$  and $\gll^s_{Q^s(R),t}=0$.
\item $\pi (\CL^s(R))= \CL^s(R/ \glsRt)$ and $\CL^s(R) = \pi^{-1} (\CL^s(R/ \glsRt))$.
    \item $Q^s(R/ \glsRt ) = Q(R/ \glsRt )$.
\end{enumerate}
\end{theorem}

{\bf Semisimplicity criterion for the ring $Q^s(R)$}. A ring is called a {\em  Goldie ring} if it a left and right Goldie ring.

\begin{theorem}\label{TY20Sep13}%\marginpar{TY20Sep13}
Let $R$ be a ring. The following statements are equivalent.
\begin{enumerate}
\item $Q^s(R)$ is a semisimple ring.
\item $R/ \glsRt $ is a semiprime  Goldie ring.
\item $Q(R/ \glsRt )$ is a semisimple ring.
\item $Q_{l, cl}(R/ \glsR )\simeq Q_{r, cl}(R/ \glsR )$ is a semisimple ring.
\end{enumerate}
If one of the equivalent conditions holds then $$Q^s(R)\simeq  Q_l(R/ \glsRt ) \simeq Q_{l,cl}(R/ \glsRt )\simeq Q_{r,cl}(R/ \glsRt ),$$  $T(R) = \pi^{-1} (\CC_{R/ \glsRt })$ and $T(R/\glsRt )= \CC_{R/\glsRt }$ where $\pi : R\ra R/ \glsRt $, $r\mapsto \br = r+\glsRt$, and $\CC_{R/ \glsRt }$ is the set of regular elements of the ring $R/ \glsRt $.
\end{theorem}

\begin{theorem}\label{TA21Sep13}%\marginpar{TA21Sep13}
Let $R=\prod_{i=1}^n R_i$ be a direct product of rings. Then
\begin{enumerate}
\item $\CL^s(R)=\prod_{i=1}^n\CL^s(R_i)$.
\item $T(R) = \prod_{i=1}^n T(R_i)$ and $ Q^s(R) \simeq \prod_{i=1}^nQ^s(R_i)$.
\end{enumerate}
\end{theorem}

\begin{theorem}\label{TA29Sep13}%\marginpar{TA29Sep13}
Let $R$ be a ring. Then $S_{l,r,\ga}(R)\subseteq T(R)$ for all $\ga \in \Ass (R)$ with $\ga \subseteq \glsRt$.
\end{theorem}

%%%%%%%%%%%%%%%%%% SECTION  5 %%%%%%%%%%%%%%%%%%%%%%%%

\section{Examples}\label{EXAMPQ}%\marginpar{EXAMPQ}

In this section, the strong left quotient ring $Q_l^s(R)$, the strong left localization radical $\gll_R^s$ and the largest strong left denominator set $T_l(R)$ are  explicitly  found for the following classes of rings: semiprime left Goldie rings (Theorem \ref{A29Aug14});  rings of $n\times n$ lower/upper triangular matrices with coefficients in a left Goldie domain (Theorem \ref{B29Aug14} and Theorem \ref{C29Aug14}); left Artinian rings (Theorem \ref{D29Aug14}), and  rings with left Artinian left quotient ring  (Theorem \ref{F29Aug14}).

{\bf Semiprime left Goldie rings}. For each semiprime left Goldie ring $R$, the theorem below describes its largest strong left quotient ring $Q_l^s(R)$, $\gll^s_R$ and $T_l(R)$.

\begin{theorem}\label{A29Aug14}%\marginpar{A29Aug14}
Let $R$ be a semiprime left Goldie ring. Then
\begin{enumerate}
\item $Q_l^s(R)= Q_{l,cl}(R)$.
\item $\gll^s_R=0$.
\item $T_l(R)=\CC_R$.
\end{enumerate}
\end{theorem}

{\it Proof}. The statements follow from Goldie's Theorem and Theorem \ref{27Sep13}.
$\Box $

%$\noindent $

{\bf Rings of lower/upper triangular matrices with coefficients in a left Goldie domain}. Let $R$ be a left Goldie domain, $D:= Q_{l,cl}(R)$ be its left quotient ring (it is a division ring), $L_n(R)$ and $U_n(D)$ be respectively the ring of lower and upper triangular matrices with coefficients in $R$ and $D$.  There are natural inclusions $R\subseteq L_n(R) \subseteq L_n(D)$ (each element $r\in R$ is identified with the diagonal matrix where all the diagonal elements are equal to $r$). Then $\CC_R=R\backslash \{ 0\}\in \Den_l(L_n(R), 0)$ with $\CC_R^{-1}L_n(R) \simeq L_n(D)$. Hence, $\CC_R\subseteq S_l(L_n(R))$ and $Q_{l,cl}(L_n(R))\simeq L_n(D)$ since $Q_{l,cl}(L_n(D))=L_n(D)$. Let $E_{ij}$ ($i,j=1, \ldots ,n$) be the matrix units. Every element $a=(a_{ij})\in L_n(R)$ is a unique sum $a= \sum_{1\leq i,j\leq n}a_{ij}E_{ij}$ where $a_{ij}\in R$.

%$\noindent $

\begin{theorem}\label{B29Aug14}%\marginpar{B29Aug14}
Let $R$ be a left Goldie domain. Then
\begin{enumerate}
\item $\maxDen_l(L_n(R))=\{ T_l(R)\}$ and $T_l(R)=\{ a=(a_{ij})\in L_n(R)\, | \, a_{11}\neq 0\}$.
\item $\gll^s_{L_n(R)}=\{ a=(a_{ij})\in L_n(R)\, | \, a_{11}=0\}$.
\item $Q_l^s(L_n(R))= Q_{l, cl}(R)$.
\end{enumerate}
\end{theorem}

{\it Proof}. Briefly, the statements follow from  Proposition  \ref{A8Dec12} and the following three facts

(i) $\maxDen_l(L_n(D))=\{ T_{E_{11}}\}$ where $T_{E_{11}}= \{ a=(a_{ij})\in L_n(D)\, | \, a_{11}\neq 0\}$, {\cite[Lemma 7.11.(2)]{Bav-LocArtRing}},

(ii) $T_{E_{11}}^{-1}L_n(D)\simeq D$ {\cite[Lemma 7.11.(2)]{Bav-LocArtRing}},

(iii) $ \ass (T_{E_{11}}) = (1-E_{11})L_n(D)= \{ a=(a_{ij})\in L_n(D)\, | \, a_{11}=0\}$, {\cite[Lemma 7.11.(3)]{Bav-LocArtRing}}.

In more detail, statement 1 follows from Proposition \ref{A8Dec12}.(3) and the statement (i).

The inclusion $T_l(R)\subseteq T_{E_{11}}$ implies the inclusion
$$ \gll^s_R = \ass (T_l(R))\subseteq \ass (T_{E_{11}})\cap L_n(R)= \ga := \{ a=(a_{ij})\in L_n(R)\, | \, a_{11}=0\}, $$by the statement (iii). Since $E_{11}\in T_l(R)$ and $E_{11}\ga =0$, we have the opposite inclusion $\gll^s_R\supseteq \ga$, i.e. $\gll^s_R=\ga $. This finishes the proof of statement 2.

Statement 3 follows from the statement (ii)  and   Proposition \ref{A8Dec12}.(4):
$$ Q_l^s(R)= T_l(R)^{-1}L_n(R)\simeq T_{E_{11}}^{-1} L_n(Q)\simeq D\simeq Q_{l,cl}(R).\;\; \Box$$

Let $U_n(R)$ be the ring of $n\times n$ upper triangular matrices with coefficients in $R$.

\begin{theorem}\label{C29Aug14}%\marginpar{C29Aug14}
Let $R$ be a left Goldie domain. Then
\begin{enumerate}
\item $\maxDen_l(U_n(R))=\{ T_l(R)\}$ and $T_l(R)=\{ a=(a_{ij})\in U_n(R)\, | \, a_{nn}\neq 0\}$.
\item $\gll^s_{U_n(R)}=\{ a=(a_{ij})\in U_n(R)\, | \, a_{nn}=0\}$.
\item $Q_l^s(U_n(R))= Q_{l,cl}(R)$.
\end{enumerate}
\end{theorem}

{\it Proof}. The theorem follows at once from Theorem \ref{B29Aug14} and the fact that the $R$-homomorphism
$$ U_n(R)\ra L_n(R), \;\; E_{ij}\mapsto E_{n+1-i, n+1-j}, $$ is a ring isomorphism. $\Box$

{\bf Left Artinian rings}. Before giving a proof of Theorem \ref{D29Aug14}, let us introduce notation and cite two results from \cite{Bav-LocArtRing}.  Let $R$ be a left Artinian ring, $\rad (R)$ be its radical, $\bR:= R/ \rad (R)=\prod_{i=1}^s\bR_i$ is  a direct  product of simple Artinian rings $\bR_i$, $\overline{1}_i$ be the identity element of the ring $\bR_i$. So, $1= \sum_{i=1}^s \overline{1}_i$ is the sum of orthogonal central idempotents of $\bR$, $1=\sum_{i=1}^s 1_i$ is a sum of orthogonal  idempotents of $R$ such that $1_i$ is a lifting of $\overline{1}_i$. For each non-empty set $I$ of $\{ 1, \ldots , s\}$, let $e_I:=\sum_{i\in I}1_i$,
$$ \CI_l':= \CI_l'(R):=\{ e_I\, | \, e_IR(1-e_I)=0\}.$$
The finite set $\CI_l'$ is a partially ordered set where $e_I\leq e_J$ if $I\subseteq J$.
\begin{proposition}\label{a16Mar14}%\marginpar{a16Mar14}
{\cite[Corollary 4.14]{Bav-LocArtRing}}
 Let $R$ be a left Artinian ring and $e:= \sum_{e'\in \min \, \CI'_l(R)}e'$. Then
\begin{enumerate}
\item $S_e:=\{1,e\}\in \Den_l(R, (1-e)R)$.
\item $\ass (S_e) = \gll_R$.
\item $e$ is the least upper bound of the set $\min \, \CI_l'(R)$ in $\CI_l'(R)$.
\end{enumerate}
\end{proposition}
The next theorem provides a description of the maximal left denominator sets of a left Artinian ring.
\begin{theorem}\label{B11Feb13}%\marginpar{B11Feb13}
{\cite[Theorem  4.10]{Bav-LocArtRing}} Let $R$ be a left Artinian ring. Then
\begin{enumerate}
\item $\maxDen_l(R)=\{ T_e\, | \, e\in \min \CI_l'(R)\} $ where $T_e=\{ u\in R\, | \, u+(1-e)R\in (R/(1-e)R)^*\}$.
\item $|\maxDen_l(R)|\leq s$  ($s$ is  the number of isomorphism classes of left simple $R$-modules).
    \item  $|\maxDen_l(R)|= s$ iff $R$ is a semisimple ring.
\end{enumerate}
\end{theorem}
The next theorem explicitly describes the triple  $T_l(R)$, $\gll^s_R$, $Q_l^s(R)$ for all left Artinian rings  $R$.
\begin{theorem}\label{D29Aug14}%\marginpar{D29Aug14}
Let $R$ be a left Artinian ring and $e= \sum_{e'\in \min \, \CI'_l(R)}e'$. Then
\begin{enumerate}
\item  $T_l(R)=\bigcap_{S\in \maxDen_l(R)}\; S= \{ u\in R\, | \, u+(1-e')R\in (R/ (1-e')R)^*$ for all $e'\in \min \, \CI_l'(R)\}$.
\item $\gll^s_R=\gll_R= (1-e)R$.
\item $Q_l^s(R)= R/ \gll_R^s\simeq R/ (1-e)R\simeq \prod_{e'\in \min \, \CI_l'(R)}\, R/ (1-e')R\simeq \prod_{S\in \maxDen_l(R)}S^{-1}R$.
\end{enumerate}
\end{theorem}

{\it Proof}.  2. By  Theorem \ref{B11Feb13}.(1), $e\in \CL_l^s(R)$. By Proposition  \ref{a16Mar14}.(1), $S_e=\{1,e\}\in \Den_l(R)$, and so $S_e\subseteq T_l(R)$, by the maximality of $T_l(R)$. Notice that $\gll^s_R\subseteq \gll_R$ (Lemma \ref{bb19Sep13}.(2)) and $\ass (S_e) = \gll_R$ (Proposition  \ref{a16Mar14}.(2)). Now,
$$ \gll_R\supseteq \gll_R^s=\ass \, T_l(R)\supseteq \ass (S_e) = (1-e)R\stackrel{{\rm Pr.\, \ref{a16Mar14}.(2)}}{=}\gll_R,$$
i.e. $\gll^s_R=\gll_R=(1-e)R$.

3. By Theorem \ref{21Sep13}.(2)  and statement 2,
$$ Q_l^s(R) = Q_l(R/\gll^s_R)= R/ \gll^s_R= R/(1-e)R.$$
The isomorphism $R/(1-e)R\simeq \prod_{e'\in \min \, CI_l'(R)}\, R/ (1-e')R$  follows from the decomposition {\cite[Equality (21)]{Bav-LocArtRing}}. It remains to notice that $R/(1-e')R\simeq T^{-1}_{e'}R$ and $\maxDen_l(R) =\{ T_{e'}\, | \, e'\in \min \, \CI_l'(R)\}$ (Theorem \ref{B11Feb13}.(1)).

1. Let $\pi : R\ra R/\gll^s_R$, $r\mapsto r+\gll^s_R$. Since the group of units $U$ of the ring $ R/\gll^s_R$ is a left denominator set of $ R/\gll^s_R\simeq T_l(R)^{-1}R$ and $\ass (T_l(R))=\gll_R^s$, the pre-image
$$\pi^{-1}(U) = \{ u\in R\, | \, u\in (1-e')R\in (R/ (1-e')R)^*\; {\rm  for \; all}\;  e'\in \min \, \CI_l'(R)\} = \bigcap_{S\in \maxDen_l(R)}\; S$$
(Theorem \ref{B11Feb13}.(1)) belongs to $\Den_l(R, \gll_R^s)$ (Lemma \ref{a21Sep13}), hence $T_l(R) = \pi^{-1}(U)$. $\Box$

{\bf Rings with left Artinian left quotient ring}. Let $A$ be a ring such that $R:= Q_{l,cl}(A)$ is a left Artinian ring. We keep the notation of the previous subsection.
\begin{theorem}\label{F29Aug14}%\marginpar{F29Aug14}
Let $A$ be a ring such that $R= Q_{l,cl}(A)$ is a left Artinian ring and $e= \sum_{e'\in \min \, \CI'_l(R)}e'$. Then
\begin{enumerate}
\item  $T_l(A)=\bigcap_{S'\in \maxDen_l(R)}\; S'$.
\item $\gll_A^s= A\cap \gll_R^s= A\cap \gll_R= A\cap (1-e)R$.
\item $Q_l^s(A)\simeq Q_l^s(R)= R/ \gll_R^s\simeq R/ (1-e)R\simeq \prod_{e'\in \min \, \CI_l'(R)}\, R/ (1-e')R\simeq \prod_{S\in \maxDen_l(R)}S^{-1}R\simeq \prod_{S'\in \maxDen_l(A)}S'^{-1}A$.
\end{enumerate}
\end{theorem}

{\it Proof}.  The theorem follows from Theorem \ref{E29Aug14}, Theorem \ref{D29Aug14} and Proposition \ref{A8Dec12}.(4). $\Box$

%$\noindent $

$\noindent $

$${\bf Acknowledgements}$$

 The work is partly supported by  the Royal Society  and EPSRC.

\small{

Department of Pure Mathematics

University of Sheffield

Hicks Building

Sheffield S3 7RH

UK

email: v.bavula@sheffield.ac.uk}

\end{document}